\documentclass[preprint,12pt,dcolumn]{elsarticle}
\usepackage{graphicx}
\usepackage{amsmath,amssymb,amsfonts}
\usepackage{amsthm}
\usepackage{color}
\usepackage{dcolumn}
\newcolumntype{d}[1]{D{.}{\cdot}{#1}}
\newcolumntype{.}{D{.}{.}{-1}}
\newcolumntype{,}{D{,}{,}{-1}}

\usepackage{hhline}
\usepackage{multirow}
\newtheorem{theorem}{Theorem}
\newtheorem{corollary}[theorem]{Corollary}

\newtheorem{lemma}[theorem]{Lemma}

\theoremstyle{remark}

\theoremstyle{definition}

\newfont{\bg}{cmr10 scaled\magstep4} 

\newcommand{\bigzerou}{\smash{\lower1.7ex\hbox{\bg 0}}}
\usepackage{bm}
\usepackage{upgreek}
\usepackage{longtable}
\usepackage{here}
\usepackage{color}
\usepackage{lineno}
\usepackage{blkarray, bigstrut}
\usepackage{mathrsfs}
\newcommand{\bbR}{\mathbb{R}}

\usepackage[%
 setpagesize=false,%
 bookmarks=true,%
 bookmarksdepth=tocdepth,%
 bookmarksnumbered=true,%
 colorlinks=true,%
 pdftitle={},%
 pdfsubject={},%
 pdfauthor={},%
 pdfkeywords={}%
]{hyperref}

\biboptions{authoryear}

\begin{document}

\begin{frontmatter}

\title{Variants of the slacks-based measure with assurance region and zeros in input--output data}
\cortext[correspondingauthor]{Corresponding author}
\author[a]{Atsushi Hori\corref{correspondingauthor}}
\ead{atsushi-hori@st.seikei.ac.jp}
\author[a]{Kazuyuki Sekitani}
\ead{kazuyuki-sekitani@st.seikei.ac.jp}

\address[a]{Seikei University, 3-3-1 Kichijoji-Kitamachi, Musashino-shi, Tokyo 180-8633, Japan}

\begin{abstract}
Incorporating an assurance region (AR) into the slacks-based measure (SBM) improves practicality; however, its efficiency measure may not have desirable properties, such as monotonicity.
We incorporate a closer target setting approach into the SBM with AR and a variant of the SBM with AR. 
We demonstrate that the efficiency measure with 
the hybrid approach has the same desirable properties as that without AR, and we also show that the
efficiency scores can be computed by solving linear programming problems. 
Our proposed approach can handle zeros in the observed input--output data without any data transformation or additional model modification.
\end{abstract}


\end{frontmatter}
\section{Introduction}
Data envelopment analysis (DEA) is a mathematical programming approach widely used to evaluate the efficiency of decision making units (DMUs) in various sectors.
Various of DEA models have been proposed since the seminal work of~\cite{charnes1978measuring}.
DEA has also been widely used to  benchmark assessed DMUs. 
\cite{charnes1978measuring} first proposed a DEA model, the Charnes, Cooper, and Rhodes (CCR) model.
The multiplier form of the CCR model, the dual problem of the CCR model, has been used to incorporate value judgments into DEA;~see \citet{thompson1986comparative} and \citet{allen1997weights} for the details. 
\par
DEA assumes that each DMU is a transformation process that consumes inputs and produces outputs. 
One advantage of the DEA is that evaluating the efficiency of DMUs does not require prior knowledge
of the transformation process from input to output. 
However, this advantage may be a weakness in real-world applications. 
Thus, incorporating process knowledge is natural in practice and typically involves constraints on input and/or output weights in the DEA model, and
these  constraints are referred to as assurance regions (AR).
\par
The efficiency score of the assessed DMU is obtained from the optimal value of the DEA model. 
 Hence,  the assessed DMU has a different efficiency score, which depends on the choice of DEA model.
\citet{fare1978measuring}, \citet{10.1007/978-3-642-52481-3_18}, \citet{russell2009axiomatic}, and \citet{SUEYOSHI2009764} provided axioms for the efficiency measurement and they proposed that 
the efficiency score should be nonnegative for all DMUs, a DMU with better activities should have a higher calculated efficiency score, the efficiency score of an efficient DMU should always be 1, etc.
We will define the desirable properties of the efficiency scores in a later section.
\par
As an alternative DEA model for CCR, \cite{tone2001slacks}  proposed   a slacks-based measure (SBM) that satisfies the desirable properties of  efficiency measures.
The original SBM DEA model was extended by~\cite{TONE2020926}, \cite{LIN2022102669}, \cite{aparicio2007closest}, \cite{ToneSBMMAX}, and \cite{KAO2022102525} to accept nonpositive-data and/or avoid a further target in benchmarking.
These extensions are effective tools for providing a practical target (benchmark) for a DMU.
However, this may lead to side effects due to the extensions. 
For instance, 
\cite{TONE2020926} and~\cite{LIN2022102669} suggested data transformation to handle zero data and multistage computation of efficiency scores.
\cite{aparicio2007closest} and \cite{KAO2022102525}
formulated closer target setting as a mix integer programing problem.
Moreover, efficiency measures using extended DEA models~developed by \cite{ToneSBMMAX}, \cite{KAO2022102525}, and \cite{aparicio2007closest} 
 no longer have the same desirable properties as the original SBM DEA model.
\par
\cite{tone2001slacks} incorporated AR into SBM (SBM-AR). 
The efficiency score of SBM-AR can be computed by solving linear programming (LP); however, it may be negative.
This potential issue is illustrated by a counterexample in a later section.
To overcome this issue, we  modify the SBM-AR DEA model and its variant using closer target setting approach. 
The modification contributes not only to practical benchmarking but also to avoiding negative efficiency scores.
Moreover, it satisfies the same properties of  efficiency measures as the original SBM DEA model.
By making minor model changes from the modified DEA model, we  can  handle zero data in an observed input--output vector. 
Additionally, the modified DEA model and its extensions can be solved by a series of linear programming problems.
\par
The remainder of this paper is organized as follows: 
Section~\ref{sec:prelim} provides methodological details of the SBM-AR DEA model and its variant, and shows by a counterexample in which the DEA models attain negative efficiency scores.
Section~\ref{sec:SBM-AR} proposes a closer target setting approach to the two SBM-AR DEA models.
This includes analyzing its properties and computational aspects.   
Section~\ref{sec:zeros.SBM-AR}  extends the two closer target setting SBM-AR DEA models 
to  handle zeros in an  input--output vector. 
Section~\ref{sec:example} presents an illustrative example that compares the proposed AR DEA models with existing AR models.
\section{Preliminary}\label{sec:prelim}
Suppose that $n$ DMUs exist, and the $j$th DMU, denoted by DMU$_j$, $j=1,\dots,n$, has $m$ inputs $\bm{x}_j=(x_{1j},\ldots.x_{mj})$ and $s$ outputs $\bm{y}_j=(y_{1j},\ldots,y_{sj})$.
We assume that $(\bm{x}_j,\bm{y}_j)\in \bbR^{m+s}_{++}$ for all $j=1,\ldots,n$. 
\par
For a given pair of matrices $P$ and $Q$ representing an assurance region (AR) of the input and output weights, respectively, the slacks-based measure with assurance region (SBM-AR) proposed by~\citet{tone2001slacks} is formulated as follows:
\begin{eqnarray}
 \min        && \frac{ 1-\frac{1}{m}\sum_{i=1}^m {d_i^-}/{x_i} }{ 1+\frac{1}{s}\sum_{r=1}^s {d^+_r}/{y_{r}}} \label{P0}\\
 \mbox{s.t.} &&  \sum_{j=1}^n \lambda_j \bm{x}_j - P \bm{\alpha} +  \bm{d}^-= \bm{x}, \label{P1}\\
  &&  \sum_{j=1}^n \lambda_j \bm{y}_j + Q \bm{\beta} - \bm{d}^+ = \bm{y}, \label{P2}\\
  && \bm{\lambda} \geq \bm{0},\   \bm{\alpha} \geq \bm{0},\  \bm{\beta} \geq \bm{0},\,  \bm{d}^- \geq \bm{0},\  \bm{d}^+ \geq \bm{0}, \label{P3}        
\end{eqnarray} 
which is a dual form of 
\begin{eqnarray}
 \max        && \bm{u}\bm{y}-\bm{v}\bm{x}+1  \label{D0}\\
 \mbox{s.t.} &&  \bm{u}\bm{y}_j - \bm{v}\bm{x}_j \leq 0, \ j=1,\ldots,n, \label{D1}\\
  &&  v_i \geq \frac{1}{mx_i}, \ i=1,\ldots,m \label{D2}\\
  && u_r \geq \frac{1+\bm{u}\bm{y}-\bm{v}\bm{x}}{sy_r},\  r=1,\ldots,s \\
  && \bm{v} P  \leq \bm{0},\   \bm{u}Q \leq \bm{0}. \label{D3}        
\end{eqnarray} 
Here, as stated  in Section 4.3 of~\citet{tone2001slacks}, the AR proposed by~\cite{thompson1986comparative} is formulated as the following ratio constraints  
\begin{align}
\underbar{$a$}_i \leq \frac{v_i}{v_1} \leq \bar{a}_i & \mbox{ for all } i=2,\ldots,m,\\
\underbar{$b$}_r \leq \frac{u_r}{u_1} \leq \bar{b}_r & \mbox{ for all } r=2,\ldots,s,
\end{align}
where $0<\underbar{$a$}_i\leq \bar{a}_i$ for all $i=2,\ldots,m$ and  $0<\underbar{$b$}_r\leq \bar{b}_r$ for all $r=2,\ldots,s$. 
Then, the corresponding $m\times (2m-2)$ matrix $P=[ p_{ij} ]$ in~\eqref{D3} is obtained  as follows:
\begin{align}
p_{ij} = \left\{
\begin{array}{rl}
\underbar{$a$}_{\frac{j+3}{2}} & (i=1 \text{ and } j \mbox{ is odd }) \\ 
-\bar{a}_{\frac{j+2}{2}} & (i=1\text{ and } j\mbox{ is even}) \\ 
-1          & (2i-3=j) \\
1           & (2i-2=j) \\
0           & \mbox{otherwise}
\end{array}
\right. \label{MatP}
\end{align}
and the corresponding $s\times 2(s-2)$ matrix $Q=[q_{rj}]$ in~\eqref{D3} is obtained as follows: 
\begin{align}
q_{rj} = \left\{
\begin{array}{rl}
\underbar{$b$}_{\frac{j+3}{2}} & (r=1\text{ and } j \mbox{ is odd}) \\ 
-\bar{b}_{\frac{j+2}{2}} & (r=1 \text{ and } j \mbox{ is even}) \\ 
-1          & (2r-3=j) \\
1           & (2r-2=j) \\
0           & \mbox{otherwise}
\end{array}.
\right.\label{MatQ}
\end{align}
\par
\citet{BRWZ} and  \citet{portela2003finding} proposed a variant  DEA model of SBM as follows:
\begin{eqnarray}
 \min        && \left(  1-\frac{1}{m}\sum_{i=1}^m \frac{d_i^-}{x_i} \right) \times\frac{1}{s} \left( \sum_{r=1}^s \frac{1}{ 1+{d^+_r}/{y_{r}} } \right) \label{P5}\\
 \mbox{s.t.} &&  \sum_{j=1}^n \lambda_j \bm{x}_j +  \bm{d}^-= \bm{x},\label{P6} \\
  &&  \sum_{j=1}^n \lambda_j \bm{y}_j  - \bm{d}^+ = \bm{y},\label{P7}\\
  && \bm{\lambda} \geq \bm{0},\   \bm{d}^- \geq \bm{0},\  \bm{d}^+ \geq \bm{0}.        \label{P8}
\end{eqnarray} 
This is known as the BRWZ measure DEA model.
The DEA model of  BRWZ measure with an assurance region  (BRWZ measure-AR) is formulated as follows:
\begin{eqnarray}
 \min        &&\left\{ \left.  \left(  1-\frac{1}{m}\sum_{i=1}^m \frac{d_i^-}{x_i} \right) \times\frac{1}{s} \left( \sum_{r=1}^s \frac{1}{ 1+{d^+_r}/{y_{r}} } \right)\,\right| 
  \eqref{P1}, \eqref{P2}, \eqref{P3}  \right\}. \label{BRWZ-AR}
\end{eqnarray}
Let
\begin{equation}
T := \left\{\, (\bm{x},\bm{y})\, \left|
\begin{array}{l}
\sum_{j=1}^n \lambda_j \bm{x}_j - P \bm{\alpha} \leq \bm{x},\\
\sum_{j=1}^n \lambda_j \bm{y}_j + Q \bm{\beta} \geq \bm{y},\\
\bm{\lambda} \geq \bm{0}, \ \bm{\alpha} \geq \bm{0}, \ \bm{\beta} \geq \bm{0} 
\end{array}
\right\}\right.\label{PPS}
\end{equation}
and 
\begin{equation}
\partial^s(T):= \left\{\, (\bm{x},\bm{y})\in T \, \left|
\begin{array}{l}  (\bm{x},-\bm{y}) \geq  (\bm{x}',-\bm{y}'),  \\
(\bm{x},-\bm{y}) \not=  (\bm{x}',-\bm{y}') 
\end{array}
\implies   (\bm{x}',\bm{y}')\notin T 
\right\}. \right.\label{partialS}
\end{equation}
Under the presence of production trade-off, 
$T$ is called the production possibility set (PPS), and $\partial^s(T)$ is called the strongly efficient frontier of $T$.
Here, according to \citet{razipour2020finding}, the presence of production trade-off assumes 
that $(\bm{x}-P\bm{\alpha},\bm{y}+Q\bm{\beta})\in T$ for any $\bm{\alpha}\geq \bm{0}$ and any $\bm{\beta}\geq \bm{0}$ if $(\bm{x},\bm{y})\in T$.
\par
The efficiency measure is a mapping $F:T\cap \bbR^{m+s}_{++} \to [0,1]$.
We consider three types of axioms as suitable efficiency measures defined on the full space of inputs and outputs.
Two of the axioms, indication and monotonicity, are  obvious extensions of the axioms proposed by~\citet{fare1978measuring} for input-based measures of efficiency.
We consider a weaker version of monotonicity as well.
\begin{description}
\item{Indication of Efficiency (I): } For all $(\bm{x},\bm{y})\in T\cap \bbR^{m+s}_{++}$, $F(\bm{x},\bm{y})=1$ if and only if $(\bm{x},\bm{y})\in \partial^s(T)$. 
\item{Strong monotonicity (SM): } For all pairs  $(\bm{x},\bm{y})\in T\cap \bbR^{m+s}_{++}$ and  $(\bm{x}',\bm{y}')\in T\cap \bbR^{m+s}_{++}$ satisfying  $(\bm{x},-\bm{y}) \leq  (\bm{x}',-\bm{y}')$ and 
$(\bm{x},-\bm{y}) \not=  (\bm{x}',-\bm{y}')$, $F(\bm{x},\bm{y}) > F(\bm{x}',\bm{y}').$
\item{Weak monotonicity (WM): } For all pairs  $(\bm{x},\bm{y})\in T\cap \bbR^{m+s}_{++}$ and  $(\bm{x}',\bm{y}')\in T\cap \bbR^{m+s}_{++}$ satisfying  $(\bm{x},-\bm{y}) \leq  (\bm{x}',-\bm{y}')$, $F (\bm{x},\bm{y}) \geq F(\bm{x}',\bm{y}').$   
\end{description}
\par
Let 
\begin{align}
g_-(\bm{d}^-) &:= 1-\frac{1}{m}\sum_{i=1}^m \frac{d^-_i}{x_i}, \\ 
g_+(\bm{d}^+) &:= \frac{s}{\displaystyle\sum_{r=1}^s\left( 1+\frac{d^+_r}{y_r}\right)}\\
\intertext{ and }
h_+(\bm{d}^+) &:= \frac{1}{s}\sum_{r=1}^s  \frac{1}{1+{d^+_r}/{y_r}}.
\end{align}
Then, the objective function of SBM-AR~\eqref{P0}--\eqref{P3} is denoted by the product of the two  functions, $g_-(\bm{d}^-)\times g_+(\bm{d}^+)$,
and that of BRWZ measure~\eqref{P5}--\eqref{P8} and BRWZ measure-AR~\eqref{BRWZ-AR} are also denoted by the product $g_-(\bm{d}^-)\times h_+(\bm{d}^+)$.
Each DEA model of  SBM-AR and BRWZ measure-AR  has a dimensionless objective function; hence, the two measures satisfy an axiom of invariance with respect to the units of measurement. 
This is an extension of the commensurability property proposed by~\citet{10.1007/978-3-642-52481-3_18} for an input-based efficiency measure. 
\par
\begin{theorem}\label{Theorem1}
Suppose  $(\bm{x},\bm{y})\in T \cap \bbR^{m+s}_{++}$ and let  $(\bm{\lambda}^*,\bm{\alpha}^*,\bm{\beta}^*,\bm{d}^{-*},\bm{d}^{+*})$ be an optimal solution to SBM-AR DEA model~\eqref{P0}--\eqref{P3}. If 
$\bm{d}^{-*}<\bm{x}$, then 
$(\bm{x}-\bm{d}^{-*},\bm{y}+\bm{d}^{+*})\in \partial^s(T)$.
\end{theorem} 
See~\ref{app:th1th2} for the proof of Theorem~\ref{Theorem1}.
\begin{theorem}\label{Theorem2}
Suppose  $(\bm{x},\bm{y})\in T \cap \bbR^{m+s}_{++}$ and let  $(\bm{\lambda}^*,\bm{\alpha}^*,\bm{\beta}^*,\bm{d}^{-*},\bm{d}^{+*})$ be an optimal solution to BRWZ measure-AR DEA model~\eqref{BRWZ-AR}.
If $\bm{d}^{-*}<\bm{x}$, then 
$(\bm{x}-\bm{d}^{-*},\bm{y}+\bm{d}^{+*})\in \partial^s(T)$.
\end{theorem} 
See~\ref{app:th1th2} for the proof of Theorem~\ref{Theorem2}.
\bigskip
\par
Following~\citet{razipour2020finding}, $(\bm{x}-\bm{d}^{-*},\bm{y}+\bm{d}^{+*})$ of Theorems~\ref{Theorem1} and~\ref{Theorem2} is called a projection point of $(\bm{x},\bm{y})$.
From $(\bm{d}^{-*},\bm{d}^{+*})\in \bbR^{m+s}_+$, the projection point dominates  $(\bm{x},\bm{y})$. 
Theorems~\ref{Theorem1} and~\ref{Theorem2} imply that optimal values of SBM-AR and BRWZ measure-AR satisfy (I).
\par
We illustrate SBM-AR DEA model~\eqref{P0}--\eqref{P3} and BRWZ measure-AR DEA model~\eqref{BRWZ-AR} with a simple dataset comprising five DMUs with two inputs and two outputs, 
 which is the dataset of~\citet{tone2001slacks} obtained by replacing  the input--output vector $(6,3,2,3)$ of DMU$_B$ with $(6,20,1,1)$. 
The weight restrictions  (58) in~\citet{tone2001slacks}  are equal to $1\leq {v_2}/{v_1} \leq 2$ and  $1\leq {u_2}/{u_1} \leq 2$, represented  by \eqref{D3}  with 
 \begin{align}
P = Q= \left[ \begin{array}{rr}
1 & -2 \\
-1   & 1
\end{array}
\right].\label{eqPQ}
\end{align}
As summarized in Table 3, \citet{tone2001slacks}  reported that three input--output vectors, $(6,3,2,2)$, DMU$_A$  and DMU$_D$ are inefficient on  the PPS without production trade-offs.
 Since $(-6,-3,2,2)\geq (-10,-10,1,1)$, the input--output vector  $(6,20,1,1)$ of DMU$_B$ also is inefficient in $T$. 
Similarly, both DMU$_A$ and DMU$_D$ are inefficient in $T$.
The dual form of SBM-AR for DMU$_E$ has an optimal value $1$ and optimal solution $(\bm{v}^*;\bm{u}^*)=(0.753,0.7535;	0.52,1)$.
Hence, DMU$_E$ is efficient in $T$.
\par 
DMU$_B$ has the same optimal solution for both SBM-AR and BRWZ measure-AR models:
\begin{align}
&\lambda^*_A=\lambda^*_B=\lambda^*_D=0,\, \lambda^*_C=\frac{3}{22},\,  \lambda^*_E=\frac{2}{11},\, \bm{\alpha}^*=\left(19\frac{3}{22},0\right),\bm{\beta}^*=\left(0,0\right)\\ 
&d^{-*}_1=23\frac{15}{22}, d^{-*}_2=d^{+*}_1=d^{+*}_2=0.
\end{align}
The dual problem of SBM-AR for DMU$_B$ has an optimal value $-{257}/{264}$ and an optimal solution
\[
\bm{v}^*=\frac{1}{12}(1,1) \mbox{ and } \bm{u}^*=\left(\frac{1}{11}, \frac{9}{88}\right).
\]
Therefore, the SBM-AR DEA model  has a negative efficiency score of $-{257}/{264}$ for DMU$_B$.
If  $d^{+}_r=0$ for all $r=1,\ldots,s$, then  the objective function value  $g_-(\bm{d}^-)\times h_{+}(\bm{d}^+)$
of the BRWZ measure-AR  DEA model coincides with that of the SBM-AR DEA model,  $g_-(\bm{d}^-)\times g_{+}(\bm{d}^+)$.
The BRWZ measure-AR DEA model  also has  a negative efficiency score for DMU$_B$.
\begin{table}[htb]
\centering
\small
\setlength{\tabcolsep}{6pt}
\renewcommand{\arraystretch}{0.91}
\caption{Illustrative example}\label{Table1}
\begin{tabular}{rrp{12pt}rrrp{12pt}rrr}\hline\hline
D&I-O& &\multicolumn{3}{c}{SBM-AR}&&\multicolumn{3}{c}{\small BRWZ measure-AR}\\[-1pt]  \cline{4-6} \cline{8-10} 
M&DATA  &&  Score &  Diff. & Rate  &&  Score & Diff. & Rate \\  \cline{4-4} \cline{8-8} 
U      &&      &\footnotesize Projection&   &            &&\footnotesize Projection&        &       \\
\hline
A &&& 0.793   &  &  & & 0.817 &  &  \\ \cline{4-4} \cline{8-8} 
$x_1$ & 4 & & 4 & 0 & 0 & & 4 & 0 & 0 \\
$x_2$ & 3 & & 3 & 0 & 0 & & 2.643 & 0.357 & 0.119 \\
$y_1$ & 2 & & 3.042 & 1.042 & 0.521 & & 2.714 & 0.714 & 0.357 \\
$y_2$ & 3 & & 3 & 0 & 0 & & 3 & 0 & 0 \\ \hline
B &&& -0.973    &  &  & & -0.973 &  &  \\ \cline{4-4} \cline{8-8} 
$x_1$ & 6 & & -17.682 & 23.682 & 3.947 & & -17.682 & 23.682 & 3.947 \\
$x_2$ & 20 & & 20 & 0 & 0 & & 20 & 0 & 0 \\
$y_1$ & 1 & & 1 & 0 & 0 & & 1 & 0 & 0 \\
$y_2$ & 1 & & 1 & 0 & 0 & & 1 & 0 & 0 \\ \hline
C &&& 1    &  &  & & 1 &  &  \\ \cline{4-4} \cline{8-8} 
$x_1$ & 8 & & 8 & 0 & 0 & & 8 & 0 & 0 \\
$x_2$ & 1 & & 1 & 0 & 0 & & 1 & 0 & 0 \\
$y_1$ & 6 & & 6 & 0 & 0 & & 6 & 0 & 0 \\
$y_2$ & 2 & & 2 & 0 & 0 & & 2 & 0 & 0 \\ \hline
D &&& 0.667    &  &  & & 0.686 &  &  \\ \cline{4-4} \cline{8-8} 
$x_1$ & 8 & & 8 & 0 & 0 & & 8 & 0 & 0 \\
$x_2$ & 1 & & 1 & 0 & 0 & & 0.435 &0.565&0.565  \\
$y_1$ & 6 & & 6 & 0 & 0 & & 6 & 0 & 0 \\
$y_2$ & 1 & & 2 & 1 & 1 & & 1.095 & 0.095 & 0.095  \\ \hline
E &&& 1    &  &  & & 1 &  &  \\ \cline{4-4} \cline{8-8} 
$x_1$ & 2 & & 2 & 0 & 0 & & 2 & 0 & 0 \\
$x_2$ & 4 & & 4 & 0 & 0 & & 4 & 0 & 0 \\
$y_1$ & 1 & & 1 & 0 & 0 & & 1 & 0 & 0 \\
$y_2$ & 4 & & 4 & 0 & 0 & & 4 & 0 & 0 \\ \hline\hline
\end{tabular}
\end{table} 
\section{A hybrid approach to SBM and BRWZ measures with AR }\label{sec:SBM-AR}
To overcome the case in which, for some $P$ and $Q$, the optimal values of SBM-AR and BRWZ measure-AR yield a negative efficiency score,
we incorporate a closer target setting approach~\citep[e.g.,][]{portela2003finding,aparicio2007closest} into SBM-AR DEA model~\eqref{P0}--\eqref{P3} and BRWZ measure-AR DEA model~\eqref{BRWZ-AR} and 
develop two efficiency measures that satisfy (I) as follows:  
\begin{align}
&\max \left\{  g_{-}(\bm{d}^-) \times  g_+(\bm{d}^+)  \, \left| \begin{array}{l}
(\bm{x}-\bm{d}^-,\bm{y}+\bm{d}^+)\in \partial^s(T),\\
 (\bm{d}^-,\bm{d}^+)\geq (\bm{0},\bm{0})  
 \end{array}\right.  \right\} \label{maxSBM-AR} \\
\intertext{ and }
&\max \left\{ \,  g_{-}(\bm{d}^-) \times h_+(\bm{d}^+)\, \left| \begin{array}{l}
(\bm{x}-\bm{d}^-,\bm{y}+\bm{d}^+)\in \partial^s(T),\\
 (\bm{d}^-,\bm{d}^+)\geq (\bm{0},\bm{0})  
 \end{array}\right.  \right\}. \label{maxBRWZ-AR}
\end{align}
In fact, as summarized in Tabel~\ref{Table1}, two models~\eqref{maxSBM-AR} and~\eqref{maxBRWZ-AR} provide positive efficiency scores for all DMUs.
\par
The closest target setting approach is formulated as the least distance DEA models~(see \cite{aparicio2016survey}),  which includes the following measures: an input-oriented least distance inefficiency measure
\begin{align}
&D^-(\bm{x},\bm{y}):=\min \left\{  \sum_{i=1}^m\frac{ d_i^-}{x_i} \,\left|   \begin{array}{l}
(\bm{x}-\bm{d}^- ,\bm{y})\in \partial^s(T),\\ 
\bm{d}^- \geq \bm{0} 
 \end{array} 
\right. \right\}, \label{LeastInput} 
\intertext{  an output-oriented least distance inefficiency measure }
&D^+(\bm{x},\bm{y}):=\min \left\{  \,  \sum_{r=1}^s \frac{d^+_r}{y_r} \, \left| \begin{array}{l}
(\bm{x},\bm{y}+\bm{d}^+)\in \partial^s(T),\\
 \bm{d}^+\geq \bm{0}  
 \end{array}\right.
  \right\}.  \label{LeastOutput}
\end{align}
\par
The aforementioned oriented least distance DEA models~\eqref{LeastInput} and~\eqref{LeastOutput} restrict their projection points to the strongly efficient frontier $\partial^s(T)$ of $T$ while in~\cite{briec1999holder}, the projection point  restriction from $\partial^s(T)$  to   the weakly  efficient frontier  of a PPS $T$ is weakened as follows: 
\begin{equation}
\partial^w(T):= \left\{\, (\bm{x},\bm{y})\in T \, \left|
\begin{array}{l}  (\bm{x},-\bm{y}) >  (\bm{x}',-\bm{y}')   
\end{array}
\Longrightarrow   (\bm{x}',\bm{y}')\notin T 
\right\}.\right.\label{partialW}
\end{equation}
By the definition, it is trivial that $\partial^s(T)\subseteq \partial^w(T)$.
Through the proper choice of matrices  $P$ and $Q$, we have  $\partial^s(T)\cap\left( \bbR^m\times
\left( \bbR^s_{+}\setminus \{\bm{0}\}\right)\right)=\partial^w(T)\cap \left( \bbR^m \times \left( \bbR^s_{+}\setminus \{\bm{0}\}\right)\right)$ as follows:
\begin{lemma}\label{Le3.ineff}
Consider matrix $P$ with $m$ rows and  matrix $Q$ with $s$ rows.
If 
\begin{align}
&\left\{ \bm{v}^{\top}\in \bbR^{m}_+ \left|\, \bm{v} P \leq \bm{0} \right.\right\} \subseteq \bbR^{m}_{++} \cup \left\{ \bm{0} \right\} \label{Vcon} \\
\intertext{and}
&
\emptyset \not= \left\{ \bm{u}^{\top}\in \bbR^s_+  \left|\, \bm{u} Q \leq \bm{0}, \, \bm{u}\not=\bm{0}\, \right.\right\} \subseteq \bbR^{s}_{++}, \label{Ucon}
\end{align}
then 
    \begin{equation}
        \partial^s(T)\cap \left( \bbR^m \times \left( \bbR^s_{+}\setminus \{\bm{0}\}\right)\right)=\partial^w(T)\cap \left( \bbR^m \times \left( \bbR^s_{+}\setminus \{\bm{0}\}\right)\right). \label{partialW=S}
    \end{equation}
\end{lemma}
See~\ref{app:le3co7} for the proof of Lemma~\ref{Le3.ineff}.
\bigskip
\par
For matrix $P$ of~\eqref{MatP}, $\bm{v}P\leq \bm{0}$ and $\bm{v}\geq \bm{0}$ satisfy $v_1=0 \Longleftrightarrow \bm{v}=\bm{0}$ and $v_1>0 \Longleftrightarrow \bm{v}>\bm{0}$, respectively;
hence,  \eqref{Vcon} is valid for matrix $P$ of~\eqref{MatP}.   
Similarly,  for matrix $Q$ of~\eqref{MatQ}, $\bm{u}Q\leq \bm{0}$ and $\bm{u} \geq \bm{0}$ satisfy $u_1=0 \Longleftrightarrow \bm{u}=\bm{0}$ and $u_1>0 \Longleftrightarrow \bm{u}>\bm{0}$. 
Since $u_1=1$ and $u_r=\underbar{$b$}_r$ $(r=2,\ldots,s)$ satisfy $\bm{u}Q\leq \bm{0}$ and $\bm{u}\not=\bm{0}$,  \eqref{Ucon} is valid for matrix $Q$ of~\eqref{MatQ}.  
\par 
The following lemma shows that the inefficiency score provided by the least distance DEA models $D^-$ and $D^+$ over $T\cap \left( \bbR^m\times \bbR^s_{++}\right)$ are invariant to the choice of the projection point restriction  between $\partial^w(T)$ and $\partial^s(T)$, and the inefficiency scores $D^-$ and $D^+$  can be computed  by solving a series of univariate linear programming problems.
\begin{lemma}\label{Le3}
Assume~\eqref{Vcon} and~\eqref{Ucon}, then
for each $(\bm{x},\bm{y})\in T \cap \bbR^{m+s}_{++}$ 
\begin{align}
D^-(\bm{x},\bm{y}) &=
\min \left\{  \sum_{i=1}^m\frac{ d_i^-}{x_i} \,\left|   \begin{array}{l}
(\bm{x}-\bm{d}^- ,\bm{y})\in \partial^w(T),\\ 
\bm{d}^- \geq \bm{0} 
 \end{array} 
\right. \right\} \label{BriecDm} \\
&= \min_{i=1,\ldots,m} \max\left\{ \,\delta \,\left|\, (\bm{x}-\delta x_i\bm{e}_i,\bm{y})\in T \right\}\right.\label{CompDm}\\
\intertext{ and }
D^+(\bm{x},\bm{y}) &=\min \left\{  \,  \sum_{r=1}^s \frac{d^+_r}{y_r} \, \left| \begin{array}{l}
(\bm{x},\bm{y}+\bm{d}^+)\in \partial^w(T),\\
 \bm{d}^+\geq \bm{0}  
 \end{array}\right.\right\} \label{BriecDp} \\
&= \min_{r=1,\ldots,s} \max\left\{ \,\delta \,\left|\, (\bm{x},\bm{y}+\delta y_r\bm{e}_r)\in T \right\}\right., \label{CompDp}
\end{align}
where $\bm{e}_i$  denotes the $i$th   unit vector, and  $\bm{e}_r$ denotes the $r$th unit vector.
\end{lemma}
See~\ref{app:le3co7} for the proof of Lemma~\ref{Le3}.
\bigskip
\par
A projection point by DEA models~\eqref{LeastInput} and~\eqref{LeastOutput}  is known as the closest target from $(\bm{x},\bm{y})$ to $\partial^s(T)$.
Let  $F_S(\bm{x},\bm{y})$ be the optimal value of SBM-AR DEA model~\eqref{maxSBM-AR}, then we have function $F_S: T\cap \bbR^{m+s}_{++}\to [0,1]$, and $F_S$ satisfies (I) and (WM) as follows:
\begin{theorem}\label{Theorem4}
Assume~\eqref{Vcon} and~\eqref{Ucon}, then efficiency measure $F_S$ satisfies (I) and (WM) and 
\begin{align}
F_S(\bm{x},\bm{y}) = \max\left\{ 1-\frac{1}{m} D^-(\bm{x},\bm{y}), \frac{1}{1+D^+(\bm{x},\bm{y})/s} \right\}\in (0,1] \label{F_S}
\end{align}
for each $(\bm{x},\bm{y}) \in T \cap \bbR^{m+s}_{++}$.
\end{theorem}
See~\ref{app:le3co7} for the proof of Theorem~\ref{Theorem4}.
\bigskip
\par
Theorem~\ref{Theorem4} represents that the projection point provided by SBM-AR DEA model~\eqref{maxSBM-AR} is the closest target. 
Furthermore, $F_S(\bm{x},\bm{y})$ can be computed by solving $(m+s)$ univariate linear programming problems by Lemma~\ref{Le3}.
\par
Let  $F_B(\bm{x},\bm{y})$ be the optimal value for the closer target approach of BRWZ measure-AR DEA model~\eqref{maxBRWZ-AR}.
Then, we have function $F_B: T\cap \bbR^{m+s}_{++}\to [0,1]$, and $F_B$ satisfies (I) and (SM) as follows.
\begin{theorem}\label{Theorem5}
Assume~\eqref{Vcon} and~\eqref{Ucon}, then efficiency measure $F_B$ satisfies (I) and (WM), and 
\begin{equation}
F_B(\bm{x},\bm{y}) = \max\left\{ 1-\frac{1}{m} D^-(\bm{x},\bm{y}), 1-\frac{1}{s}\left( \frac{D^+(\bm{x},\bm{y})}{1+D^+(\bm{x},\bm{y})}\right) \right\}\in (0,1]
\label{FB01}
\end{equation}
for each $(\bm{x},\bm{y}) \in T \cap \bbR^{m+s}_{++}$. If $m\leq s$, then $F_B$ satisfies  (SM). 
\end{theorem}
See~\ref{app:le3co7} for the proof of Theorem~\ref{Theorem5}.
\bigskip
\par
Theorem~\ref{Theorem5} indicates that the projection point provided by BRWZ measure-AR DEA model~\eqref{maxBRWZ-AR} is the closest target.
Furthermore, from~\eqref{CompDm} and~\eqref{CompDp}, BRWZ measure-AR DEA model~\eqref{maxBRWZ-AR} 
can be solved by $(m+s)$ linear programming problems. 
\par
Two functions $g_+$ and $h_+$ are the reciprocal of the arithmetic and harmonic means, respectively. Furthermore, from~$\bm{d}^+\in \bbR^{s}_+$ and $(1+{d^+_1}/{y_1},\dots,$ $1+{d^+_s}/{y_s})\in \bbR^{s}_{++}$, $g_+(\bm{d}^+)\leq h_+(\bm{d}^+)$, where 
$d_1^+=\cdots=d_s^+$ if and only if   $g_+(\bm{d}^+)= h_+(\bm{d}^+)$.
Since a feasible solution to SBM measure-AR DEA model~\eqref{maxSBM-AR} is that of BRWZ measure-AR DEA model~\eqref{maxBRWZ-AR} and vice versa, we have $F_S(\bm{x},\bm{y})\leq F_B(\bm{x},\bm{y})$ for each $(\bm{x},\bm{y})\in T\cap \bbR^{m+s}_{++}$.
The following corollary explains the inequality relationship between the efficiency scores 
of~\eqref{maxSBM-AR} and~\eqref{maxBRWZ-AR}. 
\begin{corollary}\label{Co6}
Assume~\eqref{Vcon} and~\eqref{Ucon}. For any $(\bm{x},\bm{y})\in T\cap \bbR^{m+s}_{++}$, we have
\begin{align}
F_S(\bm{x},\bm{y})\leq F_B(\bm{x},\bm{y}).
\end{align}
Moreover,
\begin{align}
&F_S(\bm{x},\bm{y})= F_B(\bm{x},\bm{y}) \Longleftrightarrow  s=1 \mbox{ or } \frac{D^-(\bm{x},\bm{y})}{m} \leq \frac{1}{s}\left(\frac{D^+(\bm{x},\bm{y}) }{1+D^+(\bm{x},\bm{y})} \right),\label{FS.eq.FB}\\
&F_S(\bm{x},\bm{y})< F_B(\bm{x},\bm{y}) \Longleftrightarrow  s\geq 2 \mbox{ and } \frac{D^-(\bm{x},\bm{y})}{m} > \frac{1}{s}\left(\frac{D^+(\bm{x},\bm{y}) }{1+D^+(\bm{x},\bm{y})} \right)\label{FS.l.FB}.
\end{align}
\end{corollary}
See~\ref{app:le3co7} for the proof of Corollary~\ref{Co6}.
\par
\section{Zeros in observed input--output data }\label{sec:zeros.SBM-AR}
In the previous sections, we assumed that $\bm{x}_j \in \bbR^m_{++}$ and   $\bm{y}_j \in \bbR^s_{++}$  for all $j=1,\ldots,n$.  
Thus, the ratios $d_i/x_{i}$ and $d_r/y_{r}$ in the objective functions, $g_{-}(\bm{d}^-) \times g_{+}(\bm{d}^+)$ and $g_{-}(\bm{d}^-)\times h_{+}(\bm{d}^+)$, in both DEA models cannot be defined for $x_{i}=0$ and $y_{r}=0$ if $\bm{x}\in\bbR^m_+\setminus\{\bm{0}\}$ and $\bm{y}\in\bbR^s_+\setminus\{\bm{0}\}$.
Now, we relax the positiveness of the input--output data  into zeros in the data, which expands the applicability of the DEA models to real world problems.   
To this end, we introduce new variables $\delta^-_i$ and $\delta^+_r$ that satisfy $d^-_i=\delta^-_ix_{i}$ and  $d^+_r=\delta^+_ry_{r}$, respectively, and we define functions as follows:
 \begin{align}
 {g}_{\bm{x}}^{\natural}(\bm{\delta}^-) & =1-\frac{1}{m}\sum_{i\in I(\bm{x})} \delta_i^-,\\
 {g}_{\bm{y}}^{\natural}(\bm{\delta}^+) & =\frac{s}{s+\sum_{r\in I(\bm{y})} \delta_r^+}, \\
 {h}_{\bm{y}}^{\natural}(\bm{\delta}^+) & =\frac{1}{s} \sum_{r\in I(\bm{y})} \frac{1}{1+\delta_r^+}+1-\frac{1}{s}|I(\bm{y})|,   
\end{align} 
where $I(\bm{z})=\left\{\, l \, \left|\,  z_l>0 \,\right\}\right.$ for $\bm{z}\in \bbR^L_{+}$.  
If  $d^-_i=\delta^-_ix_{i}$ and  $d^+_r=\delta^+_ry_{r}$, then   
\[
g_-(\bm{d}^-)= {g}_{\bm{x}}^{\natural}(\bm{\delta}^-), \ g_+(\bm{d}^+)= {g}_{\bm{y}}^{\natural}(\bm{\delta}^+) \mbox{ and }  h_+(\bm{d}^+)= {h}_{\bm{y}}^{\natural}(\bm{\delta}^+)
\mbox{ for any $(\bm{x},\bm{y})\in \bbR^{m+s}_{++}$. }  
\]
For any $(\bm{x},\bm{y})\in T\cap \bbR^{m+s}_+\setminus\{(\bm{0},\bm{0})\}$, the max SBM-AR DEA model~\eqref{maxSBM-AR}  is reformulated as follows: $F_S^{\natural}(\bm{x},\bm{y}):=$
 \begin{align}
&\max \left\{  g_{\bm{x}}^{\natural}(\bm{\delta}^-) \times  g_{\bm{y}}^{\natural}(\bm{\delta}^+)  \, \left| \begin{array}{l}
(\bm{x}-\sum_{i\in I(\bm{x})} \delta_i^-x_i\bm{e}_i, \bm{y}+\sum_{r\in I(\bm{y})} \delta_r^+y_r\bm{e}_r)\in \partial^s(T),\\
\delta_i^-=0 \ \forall i \notin I(\bm{x}),\  \delta_r^+=0 \ \forall r \notin I(\bm{y}), \\
 (\bm{\delta}^-,\bm{\delta}^+)\geq (\bm{0},\bm{0})
 \end{array}\right.  \right\} \label{maxSBM-ARnatural} \\
\intertext{ and the max BRWZ measure-AR DEA model~\eqref{maxBRWZ-AR}  is reformulated as  $F_B^{\natural}(\bm{x},\bm{y}):=$ }
&\max \left\{ \,  g_{\bm{x}}^{\natural}(\bm{\delta}^-) \times h_{\bm{y}}^{\natural}(\bm{\delta}^+)\, 
 \left| \begin{array}{l}
(\bm{x}-\sum_{i\in I(\bm{x})} \delta_i^-x_i\bm{e}_i, \bm{y}+\sum_{r\in I(\bm{y})} \delta_r^+y_r\bm{e}_r)\in \partial^s(T),  \\
\delta_i^-=0 \ \forall i \notin I(\bm{x}),\  \delta_r^+=0 \ \forall r \notin I(\bm{y}), \\
 (\bm{\delta}^-,\bm{\delta}^+)\geq (\bm{0},\bm{0})
 \end{array}\right. \right\}.  \label{maxBRWZ-ARnatural}
\end{align}
 \par
For any $(\bm{x},\bm{y})\in T\cap \bbR^{m+s}_+\setminus\{(\bm{0},\bm{0})\}$, an input-oriented least distance inefficiency measure is reformulated as follows:
\begin{align}
D^{\natural-}(\bm{x},\bm{y}) &=\min \left\{  \sum_{i\in I(\bm{x})}\delta^-_i \,
\left| 
\begin{array}{l}
(\bm{x}-\sum_{i\in I(\bm{x})} \delta_i^-x_i\bm{e}_i, \bm{y})\in \partial^s(T),\\
\delta_i^-=0 \ \forall i \notin I(\bm{x}), \  \bm{\delta}^-\geq \bm{0}
\end{array}
\right. \right\}, \label{LeastInputS} 
\intertext{ and an output-oriented least distance inefficiency measure is reformulated as follows:}
D^{\natural+}(\bm{x},\bm{y}):&=\min \left\{  \,  \sum_{r\in I(\bm{y})} \delta^+_r \, \left| 
\begin{array}{l}
(\bm{x},\bm{y}+\sum_{r \in I(\bm{y})} \delta^+_r y_r \bm{e}_r )\in \partial^s(T),\\
\delta_r^+=0 \ \forall r \notin I(\bm{y}),\ \bm{\delta}^+\geq \bm{0}  
 \end{array}\right.
  \right\}.  \label{LeastOutputS}
\end{align}
As discussed in Lemma~\ref{Le3}, the input- and output-oriented least distance inefficiency measures $D^{\natural-}$ and $D^{\natural+}$ can be computed by solving a series of univariate linear programming problems.
\begin{lemma}\label{Lemma8}
Assume~\eqref{Vcon} and~\eqref{Ucon}. Suppose $(\bm{x},\bm{y})\in T \cap \left(\bbR^{m}_{+}\setminus\{\bm{0}\}\times \bbR^{s}_+\setminus\{\bm{0}\}\right)$, then 
\begin{align}
D^{\natural-}(\bm{x},\bm{y}) &=\min \left\{  \sum_{i\in I(\bm{x})}\delta^-_i \,
\left| 
\begin{array}{l}
(\bm{x}-\sum_{i\in I(\bm{x})} \delta_i^-x_i\bm{e}_i, \bm{y})\in \partial^w(T),\\
\delta_i^-=0 \ \forall i \notin I(\bm{x}), \  \bm{\delta}^-\geq \bm{0}
\end{array}
\right. \right\} \label{natrualDm}\\
&= \min_{i\in I(\bm{x})} \max\left\{\, \delta\, \left| (\bm{x}-\delta x_i \bm{e}_i,\bm{y})\in T \right\}\right. \label{natrualminmaxDm} \\
\intertext{ and  }
D^{\natural+}(\bm{x},\bm{y})&=\min \left\{  \,  \sum_{r\in I(\bm{y})} \delta^+_r \, \left| 
\begin{array}{l}
(\bm{x},\bm{y}+\sum_{r \in I(\bm{y})} \delta^+_r y_r \bm{e}_r )\in \partial^w(T),\\
\delta_r^+=0 \ \forall r \notin I(\bm{y}),\ \bm{\delta}^+\geq \bm{0}  
 \end{array}\right.
  \right\} \label{natrualDp} \\
&= \min_{r\in I(\bm{y})} \max\left\{\, \delta\, \left| (\bm{x},\bm{y}+\delta y_r\bm{e}_r)\in T \right\}.\right.  \label{natrualminmaxDp}
\end{align}
\end{lemma}
See~\ref{app:le8th11} for the proof of Lemma~\ref{Lemma8}.
\par

\begin{lemma}\label{Lemma9}
Assume~\eqref{Vcon} and~\eqref{Ucon}. 
Suppose $(\bm{x},\bm{y})\in T \cap \left(\bbR^{m}_{+}\setminus\{\bm{0}\}\times \bbR^{s}_+\setminus\{\bm{0}\}\right)$, then 
\begin{align}
F_S^{\natural}(\bm{x},\bm{y}) = 
\max\left\{ 1-\frac{1}{m} D^{\natural -}(\bm{x},\bm{y}), \frac{1}{1+D^{\natural +}(\bm{x},\bm{y})/s} \right\}{\in(0,1]} \label{FsD}
\intertext{and}
F_B^{\natural}(\bm{x},\bm{y}) = 
\max\left\{ 1-\frac{1}{m} D^{\natural -}(\bm{x},\bm{y}), 1-\frac{1}{s}\left( \frac{D^{\natural +}(\bm{x},\bm{y})}{1+D^{\natural +}(\bm{x},\bm{y})}\right) \right\}{\in(0,1]}.\label{FbD}
\end{align}
\end{lemma}
See~\ref{app:le8th11} for the proof of Lemma~\ref{Lemma9}.
\bigskip
\par

Efficiency measures $F^\natural_S$ and $F^\natural_B$ possess the same properties as demonstrated in Theorems~\ref{Theorem4} and~\ref{Theorem5}.
\begin{theorem}\label{Theorem10}
Assume~\eqref{Vcon} and~\eqref{Ucon}, then efficiency measures $F^{\natural}_S$ and  $F^{\natural}_B$ over 
$T \cap \left(\bbR^{m}_{+}\setminus\{\bm{0}\}\times \bbR^{s}_+\setminus\{\bm{0}\}\right)$
satisfy (I) and (WM). 
Furthermore, if $m\leq s$, then  $F^{\natural}_B$ over $T \cap \left(\bbR^{m}_{+}\setminus\{\bm{0}\}\times \bbR^{s}_+\setminus\{\bm{0}\}\right)$ satisfies  (SM).
\end{theorem}
See~\ref{app:le8th11} for the proof of Theorem~\ref{Theorem10}.
\par
\bigskip
The following theorem proves the continuity of the efficiency measure $F^{\natural}_S$ or $F^{\natural}_B$ on $(\bm{x},\bm{y})$ with $I(\bm{x})\cup I(\bm{y})\not= \emptyset$.
\begin{theorem}\label{Theorem11}
Suppose $(\bm{x},\bm{y})\in T \cap \left(\bbR^{m}_{+}\setminus\{\bm{0}\}\times \bbR^{s}_+\setminus\{\bm{0}\}\right)$. Let 
\begin{align}
\mathscr{E}(\bm{x},\bm{y}) :=\left\{\, (\bm{\epsilon}^-,\bm{\epsilon}^+) \left| 
\begin{array}{ll}
\epsilon^-_i=0 \, (i\in I(\bm{x}),& \epsilon_i^->0\, (i\notin I(\bm{x}))\\
\epsilon^+_r=0 \, (r\in I(\bm{y}),& \epsilon_r^+>0\, (r\notin I(\bm{y}))
\end{array}
\right\}. \right.
\end{align}
and let 
\begin{align}
\bm{x}(\bm{\epsilon}^-)  := \bm{x}+\bm{\epsilon}^- &\mbox{ and } \bm{y}(\bm{\epsilon}^+) := \bm{y}+\bm{\epsilon}^+
\end{align}
for $(\bm{\epsilon}^-,\bm{\epsilon}^+)\in \bbR^{m+s}$.
Assume~\eqref{Vcon} and~\eqref{Ucon}, then 
\begin{align}
D^{\natural-}(\bm{x},\bm{y}) =  \lim_{ \mathscr{E}(\bm{x},\bm{y})\ni (\bm{\epsilon}^-,\bm{\epsilon}^+)  \to (\bm{0},\bm{0})}
D^-(\bm{x}(\bm{\epsilon}^-),\bm{y}(\bm{\epsilon}^+))
\intertext{and}
D^{\natural+}(\bm{x},\bm{y}) =  \lim_{ \mathscr{E}(\bm{x},\bm{y})\ni (\bm{\epsilon}^-,\bm{\epsilon}^+) \to (\bm{0},\bm{0})}
D^+(\bm{x}(\bm{\epsilon}^-),\bm{y}(\bm{\epsilon}^+)) 
\end{align}
Moreover, 
\begin{align}
F_S^{\natural}(\bm{x},\bm{y}) =  \lim_{  \mathscr{E}(\bm{x},\bm{y})\ni  (\bm{\epsilon}^-,\bm{\epsilon}^+)
\to (\bm{0},\bm{0})}
F_S(\bm{x}(\bm{\epsilon}^-),\bm{y}(\bm{\epsilon}^+))
\label{FSnat}
\intertext{and}
F_B^{\natural}(\bm{x},\bm{y}) =  \lim_{\mathscr{E}(\bm{x},\bm{y}) \ni (\bm{\epsilon}^-,\bm{\epsilon}^+)\to (\bm{0},\bm{0})}
 F_B(\bm{x}(\bm{\epsilon}^-),\bm{y}(\bm{\epsilon}^+)).
\label{FBnat}
\end{align}
\end{theorem}
See~\ref{app:le8th11} for the proof of Theorem~\ref{Theorem11}.
\par
\bigskip
Equations \eqref{FSnat} and \eqref{FBnat} of Theorem~\ref{Theorem11} represents that $F^{\natural}_S$ or $F^{\natural}_B$ are natural extension from  $F_S$ and $F_B$, respectively.
\section{An Illustrative example}\label{sec:example}
In the previous sections, we assumed that $P$ and $Q$ satisfy~\eqref{Vcon} and~\eqref{Ucon}, respectively.  
In this section, we demonstrate the verification of~\eqref{Vcon} and~\eqref{Ucon} for  given $P$ and $Q$. 
We also compare the proposed DEA-AR models with existing DEA-AR models through an illustrative example, including a DMU with zero input--output data.
Moreover, we present a counterexample to the strong monotonicity of the max SBM-AR. 
\par 
To verify assumptions~\eqref{Vcon} and~\eqref{Ucon},  we solve four linear programming problems as follows: 
\begin{align}
\hspace{-6pt}\mbox{for each } i=1,\ldots,m, & \min\left\{ v_i \, \left|\, v_1+\cdots+v_m=1, \, \bm{v} P \leq \bm{0}, \bm{v}^{\top}\in \bbR^m_+  \right\}\right.\label{eqVC}\\
\intertext{and} 
\hspace{-6pt}\mbox{for each } r=1,\ldots,s, & \min\left\{ u_r \, \left|\, u_1+\cdots+u_s=1, \, \bm{u}Q \leq \bm{0},\bm{u}^{\top}\in \bbR^s_+\right\}, \right.\label{eqUC}
\end{align}
respectively.
Matrix $P$ satisfies~\eqref{Vcon} if and only if, for all  $i\in \{ 1,\ldots,m\}$ each minimum problem~\eqref{eqVC}  does not have an optimal value $0$.  
Matrix $Q$ satisfies~\eqref{Ucon} if and only if, for all  $r\in \{ 1,\ldots,s\}$  each minimum problem~\eqref{eqUC}  has the positive optimal value. 
\par
For the given matrices $P$ and $Q$ of~\eqref{eqPQ}, all four minimum values are positive;
hence,  $P$ and $Q$  satisfy assumptions~\eqref{Vcon} and ~\eqref{Ucon}, respectively.
Let $T$ be the PPS generated  by the dataset in Table~\ref{Table1} with $P$ and  $Q$ of~\eqref{eqPQ}.   
For the PPS $T$, all the theorems in Sections~\ref{sec:SBM-AR} and~\ref{sec:zeros.SBM-AR} are valid.   
The data are summarized in Table~\ref{Table2}, which comprises eight DMUs by adding three DMUs, DMU$_F$, DMU$_G$ and DMU$_H$ to the dataset in Table~\ref{Table1}.
The three DMUs, DMU$_F$, DMU$_G$ and DMU$_H$ are all inefficient in the PPS $T$. 
The PPS generated  by the dataset in Table~\ref{Table2} with $P$ and  $Q$ of~\eqref{eqPQ} coincides with $T$ summarized in Table~\ref{Table1}.
\par
The computational results for the max SBM-AR and max BRWZ measure-AR are summarized in Table~\ref{Table2}. 
For the three common DMUs,  DMU$_A$, DMU$_B$ and DMU$_D$,  between the dataset in Table~\ref{Table1} and that in Table~\ref{Table2}, 
we  compare the efficiency scores of the four types of DEA models with AR.
In Table~\ref{Table2}, each  common DMU  has a lower efficiency score for SBM-AR (BRWZ measure-AR) than for max SBM-AR (max BRWZ measure-AR).
In particular, both max SBM-AR and  max BRWZ measure-AR provide  positive efficiency scores, whereas both  SBM-AR and  BRWZ measure-AR provide  negative efficiency scores. 
We confirm~\eqref{F_S} of Theorem~\ref{Theorem4} and~\eqref{FB01} of Theorem~\ref{Theorem5} for DMU$_B$.
\par     
Table~\ref{Table2} summarizes that DMU$_H$  is dominated by DMU$_F$; however, the max SBM-AR provides the same efficiency score for DMU$_H$ as that for DMU$_F$.
This is a counterexample of the strong monotonicity of the max SBM-AR.
Meanwhile, from $m=s=2$, the strong monotonic efficiency measure of the max BRWZ measure-AR DEA model    
provides a higher efficiency score for DMU$_H$ than that for DMU$_F$.  
 \par
DMU$_G$ in Table~\ref{Table2} has zero input data $x_{1G}$ and zero output data $y_{2G}$. 
In Table~\ref{Table2},  we observe that the efficiency score of the max SBM-AR  is $89/160$ which coincides with that of the max BRWZ measure-AR and the projection point of the max SBM-AR  coincides 
with that of max BRWZ measure-AR.
The max BRWZ measure-AR for DMU$_{G}$  has a projection point   $(0,10,80/9,0)$ instead of $(0,9/8,1,0)$.  
Indeed, from $(0,9/8,1,0)\in \partial^s(T)$ and  $(0,10,80/9,0)=80/9\times(0,9/8,1,0)$, $(0,10,80/9,0)\in \partial^s(T)$. 
Let $\delta^-_2=(10-10)/10$, $\delta^+_1=(80/9-1)/1=71/9$ and $\delta^-_1=\delta^+_2=0$, then 
$(\bm{\delta}^-,\bm{\delta}^+)$ is a feasible solution to the natural extended max BRWZ measure-AR~\eqref{maxBRWZ-ARnatural} for DMU$_G$  and 
the objective function value is 
\[
g^{\natural}_{\bm{x}_G}(\bm{\delta}^-)\times h^{\natural}_{\bm{y}_G}(\bm{\delta}^+)=   \left( 1-\frac{1}{2}\times 0 \right) \left( \frac{1}{2}\times \frac{1}{1+71/9}+1-\frac{1}{2} \right) =\frac{89}{160}.
\] 
This implies that  $(0,10,80/9,0)$ is a projection point of max BRWZ measure-AR~\eqref{maxBRWZ-ARnatural} for DMU$_G$, and there are multiple projection points of  
max BRWZ measure-AR~\eqref{maxBRWZ-ARnatural} for DMU$_G$.
\par
The provision of similar efficiency scores and a common projection point for DMU$_G$ by max SBM-AR and max BRWZ measure-AR is not a coincidence.
Suppose that DMU$_{\bar{G}}$ has $\bm{x}_{\bar{G}}=(0,\bar{x}_2)$ and $\bm{y}_{\bar{G}}=(\bar{y}_1,0)$ satisfying $(\bm{x}_{\bar{G}},\bm{y}_{\bar{G}})\in T$, $\bar{x}_2>0$ and $\bar{y}_1>0$;
then, we can show that max SBM-AR and max BRWZ measure-AR provide the same efficiency scores and the same projection points for 
DMU$_{\bar{G}}$.
Since $\left\{ (x_2,y_1) \mid (0,x_2,y_1,0)\in \partial^s(T)\right\}=\left\{ t(9/8,1) \mid t \in \bbR_+ \right\}$,  
we have 
\begin{align}
\left\{ (x_2,y_1) \mid \partial^s(T)\right\} = \left\{ (x_2,y_1) \ \middle|\ \frac{8}{9}x_2 = y_1 \right\}  \nonumber \\
(\bar{x}_2,\bar{y}_1)\in \left\{ (x_2,y_1) \mid (0,x_2,y_1,0)\in T\right\} = \left\{ (x_2,y_1) \ \middle|\ \frac{8}{9}x_2 \geq y_1 \right\} \nonumber \\
\intertext{and}
\frac{\bar{x}_2+x'}{\bar{y}_1+8/9x'} \mbox{ is a decreasing function of $x'$ on } \left[ \frac{9}{8}\bar{y}_1, \bar{x}_2\right].   \nonumber
\end{align}
The max SBM-AR DEA model for  DMU$_{\bar{G}}$ is equivalent to 
\begin{align}
&\max\left\{ \, \left. \frac{1+x'/\bar{x}_2}{1+y'/\bar{y}_1}\, \right| \, x'\leq \bar{x}_2, \ y'\geq \bar{y}_2, (0,x',y',0)\in \partial^s(T) \right\}  \nonumber \\
&= \max\left\{ \, \left. \frac{\bar{y}_1}{\bar{x}_2} \cdot \frac{\bar{x}_2+x'}{\bar{y}_1+8/9x'}\, \right| \, x'\leq \bar{x}_2, \ \frac{8}{9}x'\geq \bar{y}_1 \right\}  \nonumber \\
&= \frac{\bar{y}_1}{\bar{x}_2} \cdot \frac{\bar{x}_2+9/8\bar{y}_1}{\bar{y}_1+\bar{y}_1} =\frac{1}{2} \left( 1+ \frac{9\bar{y}_1}{8\bar{x}_2}\right). \label{eq98}
\end{align} 
Meanwhile, the max BRWZ measure-AR  DEA model for  DMU$_{\bar{G}}$ is equivalent to
\begin{align*}
&\max\left\{ \, \left. \frac{1}{4} 
 \left( 1+x'/\bar{x}_2 \right) \left( 1+\bar{y}_1/y' \right) 
\, \right| \, x'\leq \bar{x}_2, \ y'\geq \bar{y}_2, (0,x',y',0)\in \partial^s(T) \right\}  \\
&= \max\left\{ \, \left.  \frac{1}{4} 
\frac{\left( 1+x'/\bar{x}_2 \right) \left( 8/9x'+\bar{y}_1 \right)}{8/9x'} 
\right| \, x'\leq \bar{x}_2, \ \frac{8}{9}x'\geq \bar{y}_1 \right\}  .
\end{align*}
Since
$\frac{\left( 1+x'/\bar{x}_2 \right) \left( 8/9x'+\bar{y}_1 \right)}{8/9x'}$ is a convex function of $x'$ on  $\left[ \frac{9}{8}\bar{y}_1, \bar{x}_2\right]$,
maximizing the convex function over  $\left[ \frac{9}{8}\bar{y}_1, \bar{x}_2\right]$ attains at $\frac{9}{8}\bar{y}_1$ or $\bar{x}_2$.
Then the efficiency score of the BRWZ measure-AR DEA model for DMU$_G$ is obtained as follows:
\begin{align}
&\frac{1}{4} \max\left\{ \frac{\left( 1+\bar{x}_2/\bar{x}_2 \right) \left( 8/9\bar{x}_2+\bar{y}_1 \right)}{8/9\bar{x}_2}, 
\frac{\left( 1+9/8\bar{y}_1/\bar{x}_2 \right) \left( 8/9\cdot 9/8\bar{y}_1+\bar{y}_1 \right)}{8/9\cdot 9/8\bar{y}_1 }  \right\} \nonumber \\
&=\frac{1}{4} \max\left\{  2 \left(1+\frac{9\bar{y}_1}{8\bar{x}_2} \right),    \left(1+\frac{9\bar{y}_1}{8\bar{x}_2} \right)\times 2  \right\}=\frac{1}{2} \left( 1+ \frac{9\bar{y}_1}{8\bar{x}_2}\right). \label{eq99}
\end{align} 
From~\eqref{eq98} and~\eqref{eq99},  max SBM-AR and max BRWZ measure-AR provide the same efficiency scores and the same projection points for DMU$_{\bar{G}}$ with 
$\bm{x}_{\bar{G}}=(0,\bar{x}_2)$ and $\bm{y}_{\bar{G}}=(\bar{y}_1,0)$ satisfying $(\bm{x}_{\bar{G}},\bm{y}_{\bar{G}})\in T$, $\bar{x}_2>0$ and $\bar{y}_1>0$.
Since the maximization of~\eqref{eq99} attains at both $\frac{9}{8}\bar{y}_1$ and $\bar{x}_2$,  there exist multiple projection points of max BRWZ measure-AR for DMU$_{\bar{G}}$. 
\begin{table}[htb]
\centering
\small
\setlength{\tabcolsep}{6pt}
\renewcommand{\arraystretch}{0.91}
\caption{Illustrative example}\label{Table2}
\begin{tabular}{rrp{12pt}rrrp{12pt}rrr}\hline\hline
D&I-O& &\multicolumn{3}{c}{max SBM-AR}&&\multicolumn{3}{c}{\small max BRWZ measure-AR}\\[-1pt]  \cline{4-6} \cline{8-10} 
M&DATA  &&  Score &  Diff. & Rate  &&  Score & Diff. & Rate \\  \cline{4-4} \cline{8-8} 
U      &&      &\footnotesize Projection&   &            &&\footnotesize Projection&        &       \\
\hline
A&&& 0.900&&&& 0.909&&\\	\cline{4-4} \cline{8-8} 	
$x_1$&4&& 4& 0& 0&& 4& 0& 0\\		
$x_2$&3&& 3& 0& 0&& 3& 0& 0\\		
$y_1$&2&& 2& 0& 0&& 2& 0& 0\\		
$y_2$&3&& 3.667& 0.667& 0.222&& 3.667& 0.667& 0.222\\ \hline		
B&&& 0.463	&&&& 0.526&&\\	\cline{4-4} \cline{8-8} 	
$x_1$&6&& 6& 0& 0&& 6& 0& 0\\		
$x_2$&20&&-1.500&21.500& 1.075&&20& 0& 0\\		
$y_1$&1&& 1& 0& 0&& 1& 0& 0\\		
$y_2$&1&& 1& 0& 0&&19&18&18\\		\hline
C&&& 1	&&&& 1&&\\	\cline{4-4} \cline{8-8} 	
$x_1$&8&& 8& 0& 0&& 8& 0& 0\\		
$x_2$&1&& 1& 0& 0&& 1& 0& 0\\		
$y_1$&6&& 6& 0& 0&& 6& 0& 0\\		
$y_2$&2&& 2& 0& 0&& 2& 0& 0\\		\hline
D&&& 0.930	&&&& 0.930&&\\	\cline{4-4} \cline{8-8} 	
$x_1$&8&& 6.875& 1.125& 0.141&& 6.875& 1.125& 0.141\\		
$x_2$&1&& 1& 0& 0&& 1& 0& 0\\		
$y_1$&6&& 6& 0& 0&& 6& 0& 0\\		
$y_2$&1&& 1& 0& 0&& 1& 0& 0\\		\hline
E&&& 1	&&&& 1&&\\	\cline{4-4} \cline{8-8} 	
$x_1$&2&& 2& 0& 0&& 2& 0& 0\\		
$x_2$&4&& 4& 0& 0&& 4& 0& 0\\		
$y_1$&1&& 1& 0& 0&& 1& 0& 0\\		
$y_2$&4&& 4& 0& 0&& 4& 0& 0\\		\hline
F&&& 0.500	&&&& 0.530&&\\	\cline{4-4} \cline{8-8} 	
$x_1$&3&& 3& 0& 0&& 3& 0& 0\\		
$x_2$&20&& 0&20& 1&&20& 0& 0\\		
$y_1$&1&& 1& 0& 0&& 1& 0& 0\\		
$y_2$&1&& 1& 0& 0&&16.750&15.750&15.750\\		\hline
G&&& 0.556	&&&& 0.556&&\\	\cline{4-4} \cline{8-8} 	
$x_1$&0&& 0& 0& 0&& 0& 0& 0\\		
$x_2$&10&& 1.125& 8.875& 0.887&& 1.125 & 8.875& 0.887 \\ 
$y_1$&1&& 1& 0& 0&&                  1     &  0    & 0 \\ 
$y_2$&0&& 0& 0& 0&& 0& 0& 0\\		\hline
H&&& 0.500	&&&& 0.554&&\\	\cline{4-4} \cline{8-8} 	
$x_1$&3&& 3& 0& 0&& 3& 0& 0\\		
$x_2$&10&& 0&10& 1&&10& 0& 0\\		
$y_1$&1&& 1& 0& 0&& 1& 0& 0\\		
$y_2$&1&& 1& 0& 0&& 9.250& 8.250& 8.250\\		\hline\hline
\end{tabular}
\end{table}  
\section{Conclusion}\label{sec:conclusion}
In this paper, we first showed that the SBM-AR DEA model provides the negative efficiency score through a counterexample.  
We incorporated the closer  target setting approach into the SBM-AR and BRWZ measure-AR DEA models and overcame the issue.
Moreover, we demonstrated that those efficiency measures preserve the desirable properties, 
even if the observed input--output data contain zeros, for $P$ and $Q$ satisfying~\eqref{Vcon} and~\eqref{Ucon}, respectively.
Additionally, both the  SBM-AR and BRWZ measure-AR DEA models provide the closest targets to 
the strongly efficient frontier of the PPS.
They have a computational  advantage of obtaining the  efficiency score  and the target of each DMU  by solving at most $m+s$ linear programming problems.   
\par
From Theorem~\ref{Theorem5}, we have $1-\frac{1}{m}\leq 1-\frac{1}{s}< F_B(\bm{x},\bm{y})$ for each $(\bm{x},\bm{y})\in T\cap \bbR^{m+s}_{++}$ in the case of $m\leq s$; hence,
the projection point  of the max BRWZ measure-AR with $m\leq s $ belongs to  $\partial^s(T)\cap \bbR^{m+s}_{++}$. 
Since $\partial^s(T\cap \bbR^{m+s}_+)=\partial^s(T)\cap \bbR^{m+s}_+$, we have
\[
F_B(\bm{x},\bm{y}) = 
\max \left\{ \,  g_{-}(\bm{d}^-) \times h_+(\bm{d}^+)\, \left| \begin{array}{l}
(\bm{x}-\bm{d}^-,\bm{y}+\bm{d}^+)\in \partial^s(T\cap \bbR^{m+s}_+)\\
 (\bm{d}^-,\bm{d}^+)\geq (\bm{0},\bm{0})  
 \end{array}\right.  \right\}
\]
for each $(\bm{x},\bm{y})\in T\cap \bbR^{m+s}_{++}$
if $m \leq s$. 
The BRWZ measure-AR $F_B$ with $m\leq s $ is invariant, regardless of whether a nonnegative restriction is imposed on  the PPS $T$ of~\eqref{PPS}. 


\biboptions{authoryear}
\bibliography{myreference}
\appendix
\section{Proofs of~Theorems~\ref{Theorem1} and ~\ref{Theorem2}  }\label{app:th1th2}
\subsection{Proof of~Theorem~\ref{Theorem1}}
\begin{proof}
From $\bm{d}^{-*}<\bm{x}$ and $\bm{0} \leq \bm{d}^{+*}$, $g_-(\bm{d}^{-*})>0$ and $g_+(\bm{d}^{+*})>0$. 
The product of the two functions, $g_-(\bm{d}^-) \times  g_+(\bm{d}^+)$, is a  strongly decreasing function over $\left\{(\bm{d}^-,\bm{d}^+) \left|  g_-(\bm{d}^-)>0, \bm{d}^+\geq \bm{0}\right\}\right.$.
\par
Let $(\bm{\epsilon}^-,\bm{\epsilon}^+)\in \bbR^{m}_{+} \times \bbR^{s}_{+} \setminus \left\{ \left(\bm{0},\bm{0} \right)\right\}$ arbitrarily.
Then, from $g_-(\bm{d}^{-*})\times g_+(\bm{d}^{+*})>0$,  $g_+(\bm{d}^{+*}+\bm{\epsilon}^+)>0$ and the strongly decreasing function $g_-(\bm{d}^-)\times g_+(\bm{d}^+)$ over 
$\left\{(\bm{d}^-,\bm{d}^+) \left|  g_-(\bm{d}^-)>0, \bm{d}^+\geq \bm{0}\right\}\right.$, it follows that
\begin{align}
\max\left\{0, g_-(\bm{d}^{-*}+\bm{\epsilon}^-)\times  g_+(\bm{d}^{+*}+\bm{\epsilon}^+) \right\} <  g_-(\bm{d}^{-*}) \times g_+(\bm{d}^{+*}).  
\end{align}
Hence, we have
$g_-(\bm{d}^{-*}+\bm{\epsilon}^-) \times g_+(\bm{d}^{+*}+\bm{\epsilon}^+) < g_-(\bm{d}^{-*}) \times g_+(\bm{d}^{+*})$. 
This implies that from the minimum value $ g_-(\bm{d}^{-*}) \times g_+(\bm{d}^{+*})$  in the SBM-AR DEA model~\eqref{P0}--\eqref{P3}, 
\begin{equation}
\left(\bm{x}-\bm{d}^{-*}-\bm{\epsilon}^{-}, \bm{y}+\bm{d}^{+*}+\bm{\epsilon}^+ \right)\notin T.
\end{equation}
Therefore, $(\bm{x}-\bm{d}^{-*},\bm{y}+\bm{d}^{+*})\in \partial^s(T)$.
\end{proof}
\subsection{Proof of~Theorem~\ref{Theorem2}}
\begin{proof}
From $\bm{d}^{-*}<\bm{x}$ and $\bm{0} \leq \bm{d}^{+*}$, $g_-(\bm{d}^{-*})>0$ and $h_+(\bm{d}^{+*})>0$. 
The product of the two functions, $g_-(\bm{d}^-) \times  h_+(\bm{d}^+)$, is a  strongly decreasing function over $\left\{(\bm{d}^-,\bm{d}^+) \left|  g_-(\bm{d}^-)>0, \bm{d}^+\geq \bm{0}\right\}\right.$.
Hence, we can prove this using the same technique as that in Theorem~\ref{Theorem1}.
Thus, $(\bm{x}-\bm{d}^{-*},\bm{y}+\bm{d}^{+*})\in \partial^s(T)$.
\end{proof}
\section{Proofs of~Lemmas~\ref{Le3.ineff} and~\ref{Le3}, Theorems~\ref{Theorem4} and~\ref{Theorem5}, and Corollary~\ref{Co6} }\label{app:le3co7}
\subsection{Proof of~Lemma~\ref{Le3.ineff}}
\begin{proof}
From the definitions, $\partial^s(T)\subseteq \partial^w(T)$ is trivial;
hence, $\partial^s(T)\cap  \left( \bbR^m \times \left( \bbR^s_{+}\setminus \{\bm{0}\}\right)\right) \subseteq\partial^w(T)\cap  \left( \bbR^m \times \left( \bbR^s_{+}\setminus \{\bm{0}\}\right)\right)$.
\par
Suppose that any $(\bm{x},\bm{y})\in \partial^w(T)\cap  \left( \bbR^m \times \left( \bbR^s_{+}\setminus \{\bm{0}\}\right)\right)$ and  let $\bm{1}_m$ ($\bm{1}_s$) be an $m(s)$-dimensional all ones vector;
then, we have 
\begin{align}
0&=\max
\left\{\,\phi\, \left|\,
\begin{array}{l}
\sum_{j=1}^n \lambda_j \bm{x}_j -P\bm{\alpha}+\phi \bm{1}_m \leq \bm{x},\\
\sum_{j=1}^n \lambda_j \bm{y}_j +Q\bm{\beta}-\phi \bm{1}_s \geq \bm{y},\\
\bm{\lambda}\geq \bm{0},\ \bm{\alpha}\geq \bm{0}, \ \bm{\beta} \geq \bm{0}
\end{array}
\right\}\right.\label{eq30Primal}\\
&=\min
\left\{\,\bm{v}\bm{x}-\bm{u}\bm{y}\, \left|\,
\begin{array}{l}
\bm{v} \bm{x}_j -\bm{u}\bm{y}_j \geq 0,\ j=1,\ldots,n,\\
-\bm{v} P \geq \bm{0},\ -\bm{u} Q \geq \bm{0},\\
\bm{v}\bm{1}_m + \bm{u}\bm{1}_s = 1,\\
\bm{v} \geq \bm{0}, \ \bm{u} \geq \bm{0}
\end{array}
\right\}.\right.\label{eq31Dual}
\end{align}  
Let $(\bm{v}^*,\bm{u}^*)$ be the optimal solution to~\eqref{eq31Dual}.
Then, from  $\bm{v}^*\bm{1}_m + \bm{u}^*\bm{1}_s = 1$, $\bm{v}^*\bm{x}=\bm{u}^*\bm{y}$ and $\bm{y}\in \bbR^s_{+}\setminus \left\{\bm{0}\right\}$, we have
$\bm{v}^*\not=\bm{0}$. 
Let $\bar{\bm{u}}\in\bbR^s_{++}$ such that $\bm{u} Q \leq \bm{0}$ and 
$t =\min\left\{\left. \frac{\bm{v}^*\bm{x}_j}{\bar{\bm{u}}\bm{y}_j}\right| j=1,\ldots,n \right\}$.
Subsequently, $\frac{1}{\bm{v}^*\bm{1}_m+t\bar{\bm{u}}\bm{1}_s}(\bm{v}^*,t\bar{\bm{u}})$ is a feasible solution to~\eqref{eq31Dual}.
If $\bm{u}^*=\bm{0}$,
then we have $0=\bm{v}^*\bm{x}-\bm{u}^*\bm{y}=\bm{v}^*\bm{x}>\frac{1}{\bm{v}^*\bm{1}_m+t\bar{\bm{u}}\bm{1}_s}\left( \bm{v}^*\bm{x}-t\bar{\bm{u}}\bm{y}\right)$.
This contradicts the optimal value $0$ of~\eqref{eq31Dual}; hence, we have $\bm{u}^*\not=\bm{0}$.
By assumptions~\eqref{Vcon} and~\eqref{Ucon}, $\bm{v}^*> \bm{0}$ and  $\bm{u}^*> \bm{0}$. 
Therefore, the relationship between \eqref{eq30Primal} and ~\eqref{eq31Dual} is equivalent to
\begin{align}
&0= \min
\left\{\,\bm{v}\bm{x}-\bm{u}\bm{y}\, \left|\,
\begin{array}{l}
\bm{v} \bm{x}_j -\bm{u}\bm{y}_j \geq 0,  \ j=1,\ldots,n,\\
-\bm{v} P \geq \bm{0},\ -\bm{u} Q \geq \bm{0},\\
\bm{v} \geq \bm{1}, \ \bm{u} \geq \bm{1}
\end{array}
\right\}\right.=\\
&\max
\left\{\,\sum_{i=1}^m d_i^-+\sum_{r=1}^s d_r^+ \, \left|\,
\begin{array}{l}
\sum_{j=1}^n \lambda_j \bm{x}_j -P\bm{\alpha}+ \bm{d}^- =\bm{x},\\
\sum_{j=1}^n \lambda_j \bm{y}_j +Q\bm{\beta}-\bm{d}^+ = \bm{y},\\
\bm{\lambda}\geq \bm{0}, \bm{\alpha}\geq \bm{0}, \bm{\beta} \geq \bm{0}, \bm{d}^-\geq \bm{0}, \bm{d}^+\geq \bm{0}
\end{array}
\right\}.\right. \label{eq33Primal} 
\end{align}
The maximization problem~\eqref{eq33Primal} has an optimal value $0$; hence, $(\bm{x},\bm{y})\in \partial^s(T)$. 
Therefore, we have $\partial^s(T)\cap  \left( \bbR^m \times \left( \bbR^s_{+}\setminus \{\bm{0}\}\right)\right) = \partial^w(T)\cap  \left( \bbR^m \times \left( \bbR^s_{+}\setminus \{\bm{0}\}\right)\right)$. 
\end{proof}
\subsection{Proof of~Lemma~\ref{Le3}}
\begin{proof}
Suppose that any $(\bm{x},\bm{y})\in  T \cap \bbR^{m+s}_{++}$, and let $\bm{d}^{-\#}$ be an optimal solution to the right hand side of~\eqref{BriecDm}.
Then, it follows from~\eqref{partialW=S} of Lemma~\ref{Le3.ineff} and $\bm{y}\in \bbR^s_{++}$ that 
\[
(\bm{x}-\bm{d}^{-\#},\bm{y})\in \partial^w(T)\cap  \left( \bbR^m \times \bbR^s_{+}\setminus\{\bm{0}\} \right)= 
\partial^s(T)\cap \left( \bbR^m\times \bbR^s_{+}\setminus\{\bm{0}\} \right),
\]
which also implies that $\bm{d}^{-\#}$ is a feasible solution to
\begin{align}\label{Le4.ineff.strong}
    \min\left\{
        \sum_{i=1}^m\frac{d_i^-}{x_i}\ \middle|\ 
        \begin{aligned}
            &(\bm{x}-\bm{d}^-,\bm{y})\in\partial^s(T)\cap(\bbR^m\times\bbR^s_+\setminus\{\bm{0}\}),\\
            &\bm{d}^-\ge\bm{0}
        \end{aligned}
    \right\}.
\end{align}
Since $\partial^s(T)\cap \left( \bbR^m\times \bbR^s_{+}\setminus\{\bm{0}\} \right)\subseteq  \partial^s(T)\subseteq \partial^w(T)$, we have
\begin{align*}
 \sum_{i=1}^m\frac{ d_i^{-\#}} {x_i}\geq \eqref{Le4.ineff.strong} \geq D^-(\bm{x},\bm{y}) &\geq \min\left\{  \sum_{i=1}^m\frac{ d_i^-}{x_i} \,\left|  
  \begin{array}{l}
(\bm{x}-\bm{d}^- ,\bm{y})\in \partial^w(T)\\ 
\bm{d}^- \geq \bm{0} 
 \end{array} 
\right. 
\right\}\\
&=  \sum_{i=1}^m\frac{ d_i^{-\#}}{x_i}. 
\end{align*}
Therefore, \eqref{BriecDm} holds over  $T \cap \bbR^{m+s}_{++}$. 
Similarly, we can replace $\partial^s(T)$ with $\partial^w(T)$ in~\eqref{LeastOutput}; hence,  
$$
D^+(\bm{x},\bm{y})=\min \left\{  \,  \sum_{r=1}^s \frac{d^+_r}{y_r} \, \left| \begin{array}{l}
(\bm{x},\bm{y}+\bm{d}^+)\in \partial^w(T),\\ 
 \bm{d}^+\geq \bm{0}  
 \end{array}\right.
  \right\}. 
$$
\par
From Corollary~3 of~\cite{briec1999holder}, we have
\begin{align*}
D^-(\bm{x},\bm{y})&=\min \left\{  \sum_{i=1}^m\frac{ d_i^-}{x_i} \,\middle|\,  (\bm{x}-\bm{d}^- ,\bm{y})\in \partial^w(T),  \,
\bm{d}^- \geq \bm{0} 
\right\}\\
&=\min_{i=1,\ldots,m} \max\left\{ \delta \left|  
(\bm{x}-\delta x_i \bm{e}_i ,\bm{y})\in T 
\right. \right\} \\
\intertext{and }
D^+(\bm{x},\bm{y})&=\min \left\{  \,  \sum_{r=1}^s \frac{d^+_r}{y_r} \, \middle|\, (\bm{x},\bm{y}+\bm{d}^+)\in \partial^w(T), \,
 \bm{d}^+\geq \bm{0}  \right\} \\
 &=
\min_{r=1,\ldots,s} \max\left\{  \delta \left| (\bm{x},\bm{y}+\delta y_r \bm{e}_r )\in T \right.  \right\}. 
\end{align*}
\end{proof}
\subsection{Proof of~Theorem~\ref{Theorem4}}
\begin{proof}
Suppose that any $(\bm{x},\bm{y})\in  T \cap \bbR^{m+s}_{++}$,  and let 
\begin{align}
\bar{\delta}^-_i :=  \max\left\{ \delta \mid  
(\bm{x}-\delta x_i \bm{e}_i ,\bm{y})\in T 
\right\} \mbox{ for all } i=1,\ldots,m,\label{shortage_input}\\
\bar{\delta}^+_r :=  \max\left\{ \delta \mid
(\bm{x} ,\bm{y}+\delta y_r \bm{e}_r)\in T 
\right\} \mbox{ for all } r=1,\ldots,s,\label{shortage_output}
\intertext{and}
\Delta:= \left\{\, \left( 
\sum_{i=1}^m \mu_i^-\bar{\delta}^-_i\bm{e}_i, \sum_{r=1}^s \mu_r^+\bar{\delta}^+_r\bm{e}_r
 \right) 
\ \middle|\ 
 \begin{array}{l}
\sum_{i=1}^m \mu_i^-+  \sum_{r=1}^s \mu_r^+=1,\\
\bm{\mu}^- \geq \bm{0},   \bm{\mu}^+ \geq \bm{0} 
 \end{array} 
\right\}.\label{eqDelta}
\end{align}
Note that $\bar{\delta}^-_i\geq 0$ and $\bar{\delta}^+_r\ge0$ because $\delta=0$ in \eqref{shortage_input} and \eqref{shortage_output} is obviously feasible. 
It follows from Lemma~\ref{Le3.ineff} that  
\begin{align}
&\left(\bm{x}-\bar{\delta}_i^- x_i \bm{e}_i, \bm{y}\right)\in  \partial^w(T)\cap(\bbR^m\times\bbR^s_{++}) \subseteq  \partial^s(T) \mbox{ for all } i=1,\ldots,m, \label{eq35}\\
&\left(\bm{x}, \bm{y}+\bar{\delta}_r^+ y_r \bm{e}_r\right)\in \partial^w(T)\cap(\bbR^m\times\bbR^s_{++})  \subseteq  \partial^s(T) \mbox{ for all } r=1,\ldots,s. \label{eq36}
\end{align}
\par
Now, we show \eqref{F_S}.
Let $(\bm{d}^{-*},\bm{d}^{+*})$ be the optimal solution to SBM-AR model \eqref{maxSBM-AR}.   
If there exists $\bar{\delta}_i^-=0$ or  $\bar{\delta}_r^+=0$, then it follows from~\eqref{eq35} and~\eqref{eq36} that  $(\bm{x},\bm{y})\in \partial^s(T)$:
Thus, $(\bm{d}^{-*},\bm{d}^{+*})=(\bm{0},\bm{0})$.
Since $0\leq D^-(\bm{x},\bm{y})\leq \sum_{i=1}^m \frac{d^{-*}_i}{x_i}=0$ and $0\leq D^+(\bm{x},\bm{y})\leq \sum_{r=1}^s \frac{d^{+*}_r}{y_r}= 0$, 
we have  $F_S(\bm{x},\bm{y})=1=\max\left\{ 1-\frac{1}{m}0, \frac{1}{1+(1/s)\cdot 0}\right\}$.
Otherwise, that is, $\bar{\delta}_i\neq0$ for all $i$ and $\bar{\delta}_r\neq0$ for all $r$, let  
\begin{align*}
\bar{\mu}_i^- := \frac{1}{\bar{\delta}_i^-}\frac{d^{-*}_i}{x_i} \mbox{ for all } i=1,\ldots, m, \mbox{ and }
\bar{\mu}_r^+ := \frac{1}{\bar{\delta}_r^+}\frac{d^{+*}_r}{y_r} \mbox{ for all } r=1,\ldots, s.
\end{align*}
Then, we have 
\begin{align*}
\left(\bm{x}-\sum_{i=1}^m \bar{\mu}_i^- \bar{\delta}_i^{-}x_i \bm{e}_i, \bm{y}+\sum_{r=1}^s \bar{\mu}^+_r \bar{\delta}_r^+y_r\bm{e}_r
\right)=\left(\bm{x}-\bm{d}^{-*},\bm{y}+\bm{d}^{+*}\right)\in\partial^s(T). 
\end{align*}
Let 
\begin{align}
\sigma &:= \sum_{i=1}^m \bar{\mu}_i^-+ \sum_{r=1}^s \bar{\mu}_r^+ \nonumber\\
\mu_i^-&:= \frac{\bar{\mu}_i^-}{\sigma} \qquad \mbox{ for all } i=1,\ldots,m,\label{coeff_mu_i}
\intertext{ and}
\mu_r^+&:= \frac{\bar{\mu}_r^+}{\sigma} \qquad \mbox{ for all } r=1,\ldots,s.\label{coeff_mu_r}
\end{align} 
Then, we obtain $(\bm{\mu}^-,\bm{\mu}^+)\ge(\bm{0},\bm{0})$ and $\sum_{i=1}^m \mu_i^- + \sum_{r=1}^s \mu^+_r=1$.
It follows from~\eqref{eq35} and~\eqref{eq36} that 
\begin{equation}
\left( \bm{x}-\sum_{i=1}^m {\mu}_i^- \bar{\delta}_i^{-}x_i \bm{e}_i, \bm{y}+\sum_{r=1}^s {\mu}^+_r \bar{\delta}_r^+y_r\bm{e}_r
\right) \in T. \label{eq39}
\end{equation}
Now we assume that $\sigma <1$.
From~\eqref{coeff_mu_i},~\eqref{coeff_mu_r}, and $1/\sigma(\bar{\bm{\mu}}^-,\bar{\bm{\mu}}^+)=(\bm{\mu}^-,\bm{\mu}^+)$, 
we have $(\bar{\bm{\mu}}^-,\bar{\bm{\mu}}^+)\leq (\bm{\mu}^-,\bm{\mu}^+)$ and  $(\bar{\bm{\mu}}^-,\bar{\bm{\mu}}^+)\not= (\bm{\mu}^-,\bm{\mu}^+)$; hence, 
\begin{align*}
\left(\sum_{i=1}^m \bar{\mu}_i^- \bar{\delta}_i^{-}x_i \bm{e}_i, \sum_{r=1}^s \bar{\mu}^+_r \bar{\delta}_r^+y_r\bm{e}_r
\right) \leq \left(\sum_{i=1}^m {\mu}_i^- \bar{\delta}_i^{-}x_i \bm{e}_i, \sum_{r=1}^s {\mu}^+_r \bar{\delta}_r^+y_r\bm{e}_r
\right)\\
\intertext{and}
\left(\sum_{i=1}^m \bar{\mu}_i^- \bar{\delta}_i^{-}x_i \bm{e}_i, \sum_{r=1}^s \bar{\mu}^+_r \bar{\delta}_r^+y_r\bm{e}_r
\right) \not= \left(\sum_{i=1}^m {\mu}_i^- \bar{\delta}_i^{-}x_i \bm{e}_i, \sum_{r=1}^s {\mu}^+_r \bar{\delta}_r^+y_r\bm{e}_r
\right).
\end{align*} 
This result contradicts $\left(\bm{x}-\bm{d}^{-*},\bm{y}+\bm{d}^{+*}\right)\in\partial^s(T)$.  Therefore, we have  $\sigma \geq 1$, and 
 $(\bm{d}^{-*},\bm{d}^{+*})=\left( \sum_{i=1}^m \bar{\mu}_i^- \bar{\delta}_i^{-}x_i \bm{e}_i, \sum_{r=1}^s \bar{\mu}^+_r \bar{\delta}_r^+y_r\bm{e}_r\right)$ leads to
\begin{equation}
(\bm{d}^{-*},\bm{d}^{+*}) \geq \left(\sum_{i=1}^m {\mu}_i^- \bar{\delta}_i^{-}x_i \bm{e}_i, \sum_{r=1}^s {\mu}^+_r \bar{\delta}_r^+y_r\bm{e}_r
\right).\label{eq40}
\end{equation}
Since $(1-\frac{1}{m}\sum_{i=1}^m\delta_i^-)/(1+\frac{1}{s}\sum_{r=1}^s \delta_r^+)$ is a quasi-convex function with respect to $(\bm{\delta}^-,\bm{\delta}^+)$, from~\eqref{eq40},~\eqref{eq35},~\eqref{eq36}, \eqref{CompDm}, and~\eqref{CompDp}, it follows that
that  
\begin{align*}
F_S(\bm{x},\bm{y})
&=\frac{1-\frac{1}{m}\sum_{i=1}^m d_i^{-*}/x_i }{1+\frac{1}{s}\sum_{r=1}^s d_r^{+*}/y_r}\\
                                                                   &\leq \frac{1-\frac{1}{m}\sum_{i=1}^m\mu_i^-\bar{\delta}_i^-}{1+\frac{1}{s}\sum_{r=1}^s\mu_r^+\bar{\delta}_r^+}\\
                                                                   &\leq \max\left\{  \left.  \frac{1-\frac{1}{m}\sum_{i=1}^m\delta_i^-}{1+\frac{1}{s}\sum_{r=1}^s \delta_r^+}\right| (\bm{\delta}^-,\bm{\delta}^+)\in \Delta \right\}\\
                                                                   &=\max\left\{
                                                                   \max\left\{\left. 1-\frac{1}{m}\bar{\delta}_i^- \right|\, i=1,\ldots,m\right\}, 
                                                                   \max\left\{\left. \frac{1}{1+\bar{\delta}_r^+/s} \right|\, r=1,\ldots,s\right\} 
                                                                   \right\}\\
                                                                   &=\max\left\{
                                                                   1- \frac{1}{m}\min_{i=1,\ldots,m}\bar{\delta}_i^-, 
                                                                   \frac{1}{1+{1}/{s}\cdot \min_{r=1,\ldots,s}\bar{\delta}_r^+}  
                                                                   \right\}\\
                                                                   &=\max\left\{
                                                                    1- \frac{1}{m}D^-(\bm{x},\bm{y}),  \frac{1}{1+D^+(\bm{x},\bm{y})/s}  
                                                                   \right\}\\
                                                                   &\leq F_S(\bm{x},\bm{y}).
\end{align*}
Therefore,  for each $(\bm{x},\bm{y}) \in T \cap \bbR^{m+s}_{++}$
\begin{align}
F_S(\bm{x},\bm{y}) = \max\left\{ 1-\frac{1}{m} D^-(\bm{x},\bm{y}), \frac{1}{1+D^+(\bm{x},\bm{y})/s} \right\}.\label{eqFs=max}
\end{align}
\par
We now show that $F_S$ satisfies (I) and (WM).
It follows from~\eqref{eqFs=max}   that  $F_S(\bm{x},\bm{y})=1$ if and only if $D^-(\bm{x},\bm{y})=0$ or  $D^+(\bm{x},\bm{y})=0$ and the definitions \eqref{LeastInput} and~\eqref{LeastOutput} of $D^-$ and $D^+$ imply that $(\bm{x},\bm{y})\in \partial^s(T)$ is equivalent to $D^-(\bm{x},\bm{y})=0$ or  $D^+(\bm{x},\bm{y})=0$. Therefore,  $F_S$ satisfies (I).
\par
For each $(\bm{x},\bm{y}) \in T \cap \bbR^{m+s}_{++}$, we have $1+\frac{1}{s}D^+(\bm{x},\bm{y})\geq 1>0$ and $ 1-\frac{1}{m}D^-(\bm{x},\bm{y})\leq 1$.
Therefore, from~\eqref{eqFs=max}, we obtain 
\[
0<F_S(\bm{x},\bm{y})=\max\left\{ 1-\frac{1}{m}D^-(\bm{x},\bm{y}), \frac{1}{1+D^+(\bm{x},\bm{y})/s}\right\}\leq 1.
\]
for each $(\bm{x},\bm{y}) \in T \cap \bbR^{m+s}_{++}$.
Moreover, by~\eqref{BriecDm},~\eqref{BriecDp}, and Proposition~2 of~\cite{briec1999holder}, 
both the least distance inefficiency measures $D^-$ and $D^+$ over $T \cap \bbR^{m+s}_{++}$ satisfy 
weak monotonicity; hence,  $F_S$ also satisfies (WM).  
\end{proof}
\subsection{Proof of~Theorem~\ref{Theorem5}}
\begin{proof}
We can show~\eqref{eq35} and~\eqref{eq36} by the same proof as that of Theorem~\ref{Theorem4}.
\par
Let $(\bm{d}^{-*},\bm{d}^{+*})$ be the optimal solution to BRWZ measure-AR DEA \eqref{maxBRWZ-AR}.   
If there exists $\bar{\delta}_i^-=0$ or  $\bar{\delta}_r^+=0$, then it follows from~\eqref{eq35} and~\eqref{eq36} that $(\bm{x},\bm{y})\in \partial^s(T)$.
Thus, $(\bm{d}^{-*},\bm{d}^{+*})=(\bm{0},\bm{0})$.
Since $0\leq D^-(\bm{x},\bm{y})\leq \sum_{i=1}^m \frac{d^{-*}_i}{x_i}=0$ and $0\leq D^+(\bm{x},\bm{y})\leq \sum_{r=1}^s \frac{d^{+*}_r}{y_r}= 0$, 
we have  $F_B(\bm{x},\bm{y})=1=\max\left\{ 1-\frac{1}{m}0, 1-\frac{1}{s}\frac{0}{1+0} \right\}$.
Otherwise, a similar to that of Theorem~\ref{Theorem4}, we can show that
\begin{equation}
(\bm{d}^{-*},\bm{d}^{+*}) \geq \left(\sum_{i=1}^m {\mu}_i^- \bar{\delta}_i^{-}x_i \bm{e}_i, \sum_{r=1}^s {\mu}^+_r \bar{\delta}_r^+y_r\bm{e}_r
\right),\label{eq41}
\end{equation}
where $\bar{\delta}^-_i$, $i=1,\dots,m$, and $\bar{\delta}^+_r$, $r=1,\dots,s$, are defined by \eqref{shortage_input} and \eqref{shortage_output}, respectively, and
\begin{align*}
\mu_i^- := \frac{ \displaystyle \frac{1}{\bar{\delta}^-_i} \frac{d_i^{-*}}{x_i}}{ 
\displaystyle\sum_{k=1}^m \frac{1}{\bar{\delta}^-_k} \frac{d_k^{-*}}{x_k}+
\sum_{r=1}^s \frac{1}{\bar{\delta}^+_r} \frac{d_r^{+*}}{y_r}} \mbox{ for all } i=1,\ldots,m,\\
\intertext{and}
\displaystyle\mu_r^+ := \frac{ \displaystyle \frac{1}{\bar{\delta}^+_r} \frac{d_r^{+*}}{y_r}}{ 
\displaystyle\sum_{i=1}^m \frac{1}{\bar{\delta}^+_i} \frac{d_i^{-*}}{x_i}+
\sum_{k=1}^s \frac{1}{\bar{\delta}^+_k} \frac{d_k^{+*}}{y_k} } \mbox{ for all } r=1,\ldots,s.  
\end{align*}
Since $ \left(1-\frac{1}{m}\sum_{i=1}^m\delta_i^-\right)\left(  \frac{1}{s} \sum_{r=1}^s  \frac{1}{1+\delta_r^+} \right)$ is a quasi-convex function of $(\bm{\delta}^-,\bm{\delta}^+)$, from~\eqref{eq41},~\eqref{eq35},~\eqref{eq36},~\eqref{CompDm}, and~\eqref{CompDp}, it follows that 
\begin{align*}
F_B(\bm{x},\bm{y})
&=\left(1-\frac{1}{m}\sum_{i=1}^m\frac{d_i^{-*}}{x_i} \right) \left( \frac{1}{s}\sum_{r=1}^s\frac{1}{1+{d_r^{+*}}/{y_r}}\right)\\
                                                                   &\leq \left(1-\frac{1}{m}\sum_{i=1}^m\mu_i^-\bar{\delta}_i^-\right) \left( \frac{1}{s}\sum_{r=1}^s\frac{1}{1+\mu_r^+\bar{\delta}_r^+}\right)\\
                                                                   &\leq \max\left\{ \left. \left(1-\frac{1}{m}\sum_{i=1}^m\delta_i^-\right)\left( \frac{1}{s}\sum_{r=1}^s \frac{1}{1+\delta_r^+}\right)\right|  (\bm{\delta}^-,\bm{\delta}^+)\in \Delta \right\}\\
                                                                   &=\max\left\{
                                                                   \max\left\{\left. 1-\frac{1}{m}\bar{\delta}_i^- \right|\, i=1,\ldots,m\right\}, 
                                                                   \max\left\{\left. 1-\frac{1}{s}\frac{\bar{\delta}^+_r}{1+\bar{\delta}_r^+} \right|\, r=1,\ldots,s\right\} 
                                                                   \right\}\\
                                                                   &=\max\left\{
                                                                   1- \frac{1}{m}\min_{i=1,\ldots,m}\bar{\delta}_i^-, 
                                                                   1-\frac{1}{s}\min_{r=1,\ldots,s} \frac{\bar{\delta}_r^+}{1+ \bar{\delta}_r^+}  
                                                                    \right\}\\
                                                                   &=\max\left\{
                                                                    1- \frac{1}{m}D^-(\bm{x},\bm{y}),  1-\frac{1}{s}\left(\frac{D^+(\bm{x},\bm{y})}{1+D^+(\bm{x},\bm{y})} \right) 
                                                                   \right\}\\
                                                                   &\leq F_B(\bm{x},\bm{y}).
\end{align*}
Therefore,  for each $(\bm{x},\bm{y}) \in T \cap \bbR^{m+s}_{++}$
\begin{align}
F_B(\bm{x},\bm{y}) = \max\left\{ 1-\frac{1}{m} D^-(\bm{x},\bm{y}), 1-\frac{1}{s}\left(\frac{D^+(\bm{x},\bm{y})}{1+D^+(\bm{x},\bm{y})}\right) \right\}.\label{eqFB=max}
\end{align}
It follows from~\eqref{eqFB=max}   that  $F_B(\bm{x},\bm{y})=1$ if and only if $D^-(\bm{x},\bm{y})=0$ or  $D^+(\bm{x},\bm{y})=0$, and the definitions~\eqref{LeastInput} and~\eqref{LeastOutput} of $D^-$ and $D^+$ imply that $(\bm{x},\bm{y})\in \partial^s(T)$ is equivalent to $D^-(\bm{x},\bm{y})=0$ or  $D^+(\bm{x},\bm{y})=0$. Therefore,  $F_B$ satisfies (I).
\par
Next, we prove that $F_B$ satisfies (WM).
For each $(\bm{x},\bm{y}) \in T \cap \bbR^{m+s}_{++}$ we have
$1\geq 1-\frac{1}{s}\left(\frac{D^+(\bm{x},\bm{y})}{1+D^+(\bm{x},\bm{y})}\right)>1-\frac{1}{s}\geq 0$ and $1\geq 1-\frac{D^-(\bm{x},\bm{y})}{m}$ and hence, 
$0\leq 1- \frac{1}{s}< \max\left\{ 1-\frac{D^-(\bm{x},\bm{y})}{m}, 1-\frac{1}{s}\left( \frac{D^+(\bm{x},\bm{y})}{1+D^+(\bm{x},\bm{y})}\right) \right\}\leq 1$.
It follows from~\eqref{eqFB=max} that $0\leq 1-\frac{1}{s}<F_B(\bm{x},\bm{y})\leq 1$ for each $(\bm{x},\bm{y}) \in T \cap \bbR^{m+s}_{++}$.
Moreover, it follows from~\eqref{BriecDm} and \eqref{BriecDp} and Proposition~2 of~\cite{briec1999holder}, 
both the least distance inefficiency measures $D^-$ and $D^+$ over $T \cap \bbR^{m+s}_{++}$ satisfy  weak monotonicity and hence,  $F_B$ also satisfies (WM).  
\par
We demonstrate that $F_B$ has (SM) under $m\le s$.
For each $(\bm{x}',\bm{y}')\in T\cap \bbR^{m+s}_{++}$ satisfying $(\bm{x},-\bm{y})\leq (\bm{x}',-\bm{y}')$ and  $(\bm{x},-\bm{y})\not= (\bm{x}',-\bm{y}')$, let
\begin{align}
\bar{\delta}'^-_i :=  \max\left\{ \delta \mid  
(\bm{x}'-\delta x'_i \bm{e}_i ,\bm{y}')\in T 
\right\} \mbox{ for all } i=1,\ldots,m, \label{eq42}\\
\bar{\delta}'^+_r :=  \max\left\{ \delta \mid  
(\bm{x}' ,\bm{y}'+\delta y'_r \bm{e}_r)\in T 
\right\} \mbox{ for all } r=1,\ldots,s.\label{eq43}
\end{align}
Then, there exists an index $r'\in \{ 1,\ldots,s\}$ such that $\bar{\delta}'^{+}_{r'}=D^+(\bm{x}',\bm{y}')$.
Note that $\bar{\delta}'^-_i\geq 0$ and $\bar{\delta}'^+_r\ge0$ since both \eqref{eq42} and \eqref{eq43} have a feasible solution $\delta=0$.
It follows from Lemma~\ref{Le3.ineff} that
$(\bm{x}' ,\bm{y}'+\bar{\delta}'^+_{r'} y'_{r'} \bm{e}_{r'})\in \partial^w(T)\cap (\bbR^m\times\bbR^s_{++})=\partial^s(T)\cap(\bbR^m\times\bbR^s_{++}) \subseteq \partial^s(T)$.
This implies that  
\[
(\bm{x}',\bm{y}'+\bar{\delta}'^+_{r'}y'_{r'}\bm{e}_{r'})\in \partial^s(T) \mbox{ and  } (\bm{x},\bm{y}+\bar{\delta}^+_{r'}y_{r'}\bm{e}_{r'})\in\partial^s(T),
\]
where the second set membership holds from~\eqref{eq35} and~\eqref{eq36}.
From  $(\bm{x},-\bm{y})\leq (\bm{x}',-\bm{y}')$ and  $(\bm{x},-\bm{y})\not= (\bm{x}',-\bm{y}')$, 
$y'_{r'}(1+\bar{\delta}'^+_{r'})=y'_{r'}+\bar{\delta}'^+_{r'}y'_{r'}>y_{r'}+\bar{\delta}^+_{r'}y_{r}'= y_{r'}(1+\bar{\delta}^+_{r'})$.
This implies from $y'_{r'}\leq y_{r'}$ that $1\geq \frac{y'_{r'}}{y_{r'}}> \frac{1+\bar{\delta}^+_{r'}}{1+\bar{\delta}'^+_{r'}}$.
Therefore, we have $
\bar{\delta}^+_{r'} < \bar{\delta}'^+_{r'}$ and  $D^+(\bm{x},\bm{y})\leq \bar{\delta}^+_{r'}< \bar{\delta}'^+_{r'}=D^+(\bm{x}',\bm{y}').
$
\par
If $F_B(\bm{x}',\bm{y}')=1-\frac{1}{s}\left[{D^+(\bm{x}',\bm{y}')}/{(1+D^+(\bm{x}',\bm{y}'))}\right]$, then
\[
F_B(\bm{x}',\bm{y}')=1-\frac{1}{s}\left(\frac{D^+(\bm{x}',\bm{y}')}{1+D^+(\bm{x}',\bm{y}')} \right)<1-\frac{1}{s}\left(\frac{D^+(\bm{x},\bm{y})}{1+D^+(\bm{x},\bm{y})}\right)\leq F_B(\bm{x},\bm{y}).
\]
Suppose that $F_B(\bm{x}',\bm{y}')=1-\frac{1}{m}D^-(\bm{x}',\bm{y}')$; then, there exists an index $i'\in \{1,\ldots,m\}$ such that 
$\bar{\delta}'^{-}_{i'}=D^-(\bm{x}',\bm{y}')$, and we have
$1-\frac{1}{s} < F_B(\bm{x}',\bm{y}') = 1-\frac{1}{m} \bar{\delta}'^{-}_{i'}$. 
This implies that, from $m\leq s$, 
\begin{equation}
\bar{\delta}'^{-}_{i'}<1. \label{eq44}
\end{equation}
Since $(\bm{x}'-\bar{\delta}'^{-}_{i'} x'_{i'} \bm{e}_{i'} ,\bm{y}')\in \partial^s(T)$ and $(\bm{x}-\bar{\delta}^{-}_{i'} x_{i'} \bm{e}_{i'} ,\bm{y})\in \partial^s(T)$, we have
$(1-\bar{\delta}'^-_{i'})x'_{i'}<(1-\bar{\delta}^-_{i'})x_{i'}$.
Thus, from~\eqref{eq44}, $1\leq {x'_{i'}}/{x_{i'}}< ({1-\bar{\delta}^-_{i'}})/({1-\bar{\delta}'^-_{i'}})$; hence,
$\bar{\delta}^-_{i'}<\bar{\delta}'^-_{i'}$. Therefore, we have
\[
F_B(\bm{x}',\bm{y}')=1-\frac{1}{m} \bar{\delta}'^-_{i'} < 1-\frac{1}{m} \bar{\delta}^-_{i'}\leq  1-\frac{1}{m}D^-(\bm{x},\bm{y})\leq F_B(\bm{x},\bm{y}).
\]  
Efficiency measure $F_B$ satisfies (SM) if $m\leq s$.
\end{proof}
\subsection{Proof of~Corollary~\ref{Co6}}
\begin{proof}
It suffices to show \eqref{FS.eq.FB} and \eqref{FS.l.FB}.
For each $(\bm{x},\bm{y})\in T \cap \bbR^{m+s}_{++}$, it follows from $s\geq 1$ and $D^+(\bm{x},\bm{y})\geq 0$ that 
\begin{align}
1-\frac{1}{s}\left(\frac{D^+(\bm{x},\bm{y}) }{1+D^+(\bm{x},\bm{y})} \right)&-\frac{1}{1+\frac{1}{s}D^+(\bm{x},\bm{y})} =
\frac{ {D^+(\bm{x},\bm{y})}/{s} }{1+\frac{1}{s}D^+(\bm{x},\bm{y})}-\frac{ {D^+(\bm{x},\bm{y})}/{s} }{1+D^+(\bm{x},\bm{y})}\nonumber\\
&=\frac{D^+(\bm{x},\bm{y})}{s}\left(\frac{1}{1+D^+(\bm{x},\bm{y})/s}-\frac{1}{1+D^+(\bm{x},\bm{y})}  \right)\nonumber\\
&=\frac{(D^+(\bm{x},\bm{y}))^2}{s}\left(\frac{1-{1}/{s}}{ \left(1+{D^+(\bm{x},\bm{y})}/{s}\right)\left(1+D^+(\bm{x},\bm{y})\right)} \right)\nonumber\\
&\geq 0,\label{ieq:D+.frac}
\end{align}
where the equality in the last inequality equation holds if and only if $s=1$ or $D^+(\bm{x},\bm{y})= 0$. Since $D^+(\bm{x},\bm{y})= 0$ is equivalent to  $D^-(\bm{x},\bm{y})=0$, 
${D^-(\bm{x},\bm{y})}/{m} > \frac{1}{s}[{D^+(\bm{x},\bm{y})}/({1+D^+(\bm{x},\bm{y})})]$ implies both $D^+(\bm{x},\bm{y})>0$ and $D^-(\bm{x},\bm{y})>0$, which also implies
\begin{equation}
1-\frac{1}{s}\left(\frac{D^+(\bm{x},\bm{y}) }{1+D^+(\bm{x},\bm{y})} \right)>\frac{1}{1+D^+(\bm{x},\bm{y})/s}. \label{eq48}
\end{equation}
\par
If   ${D^-(\bm{x},\bm{y})}/{m} \leq \frac{1}{s}[{D^+(\bm{x},\bm{y})}/({1+D^+(\bm{x},\bm{y})})]$, then \eqref{ieq:D+.frac} leads that 
\[
1-\frac{D^-(\bm{x},\bm{y})}{m}\geq 1-\frac{1}{s}\left(\frac{D^+(\bm{x},\bm{y}) }{1+D^+(\bm{x},\bm{y})} \right)\geq \frac{1}{1+D^+(\bm{x},\bm{y})/s}.
\]
This implies from Theorems~\ref{Theorem4} and~\ref{Theorem5}  that  
\begin{align*}
F_B(\bm{x},\bm{y})&=
\max\left\{1-\frac{D^-(\bm{x},\bm{y})}{m},1-\frac{1}{s}\left(\frac{D^+(\bm{x},\bm{y}) }{1+D^+(\bm{x},\bm{y})} \right) \right\}\\
&= 1-\frac{D^-(\bm{x},\bm{y})}{m} \\
&= \max\left\{
1-\frac{D^-(\bm{x},\bm{y})}{m} , \frac{1}{1+D^+(\bm{x},\bm{y})/s} \right\} =F_S(\bm{x},\bm{y}). 
\end{align*}
\par
If $s=1$, then $1-\frac{1}{s}[{D^+(\bm{x},\bm{y}) }/({1+D^+(\bm{x},\bm{y})})]=\frac{1}{1+D^+(\bm{x},\bm{y})/s}$, which implies from Theorems~\ref{Theorem4} and~\ref{Theorem5} that
 $F_B(\bm{x},\bm{y})=F_S(\bm{x},\bm{y})$.
\par
If ${D^-(\bm{x},\bm{y})}/{m} > \frac{1}{s}[{D^+(\bm{x},\bm{y})}/({1+D^+(\bm{x},\bm{y})})]$ and $s\geq 2$, then Theorems~\ref{Theorem4}, \ref{Theorem5}, and~\eqref{eq48} show that
\begin{align*}
F_B(\bm{x},\bm{y})&=
\max\left\{1-\frac{D^-(\bm{x},\bm{y})}{m},1-\frac{1}{s}\left(\frac{D^+(\bm{x},\bm{y}) }{1+D^+(\bm{x},\bm{y})} \right) \right\}\\
&=1-\frac{1}{s}\left(\frac{D^+(\bm{x},\bm{y}) }{1+D^+(\bm{x},\bm{y})} \right) \\
&>\max\left\{
1-\frac{D^-(\bm{x},\bm{y})}{m} , \frac{1}{1+D^+(\bm{x},\bm{y})/s} \right\} \\
&=F_S(\bm{x},\bm{y}). 
\end{align*}
\end{proof}
\section{Proofs of~Lemma~\ref{Lemma8} and~\ref{Lemma9}, Theorems~\ref{Theorem10} and~\ref{Theorem11} }\label{app:le8th11}
\subsection{Proof of~Lemma~\ref{Lemma8}}
\begin{proof}
Suppose any $(\bm{x},\bm{y})\in  T \cap \left( \bbR^{m}_{+}\times \bbR^{s}_+\setminus\{\bm{0}\}\right)$ and 
let $\bm{\delta}^{-\#}$ be an optimal solution to the right-hand side of~\eqref{natrualDm}.
Then, it follows from~\eqref{partialW=S} of Lemma~\ref{Le3.ineff} and $\bm{y}\in \bbR^s_{+}\setminus\{\bm{0}\}$ that 
\[
(\bm{x}-\sum_{i\in I(\bm{x})} \delta^{-\#}_ix_i\bm{e}_i,\bm{y})\in \partial^w(T)\cap  \left( \bbR^m \times \bbR^s_{+}\setminus\{\bm{0}\} \right)= 
\partial^s(T)\cap \left( \bbR^m\times \bbR^s_{+}\setminus\{\bm{0}\} \right) ,
\]
which implies that $\bm{\delta}^{-\#}$ is feasible for
\begin{align}\label{Le8.ineff.strong}
    \min\left\{\sum_{i\in I(\bm{x})}\delta_i^-\ \middle|\ 
    \begin{aligned}
        &\textstyle(\bm{x}-\sum_{i\in I(\bm{x})}\delta_i^-x_i\bm{e}_i,\bm{y})\in\partial^s(T)\cap \left( \bbR^m\times \bbR^s_{+}\setminus\{\bm{0}\} \right) ,\\
        &\delta_i^-=0\ \forall i\notin I(\bm{x}),\ \bm{\delta}^-\ge\bm{0}
    \end{aligned}
    \right\}.
\end{align}
Since $\partial^s(T)\cap \left( \bbR^m\times \bbR^s_{+}\setminus\{\bm{0}\} \right)\subseteq  \partial^s(T)\subseteq \partial^w(T)$, we have
\begin{align*}
&\sum_{i\in I(\bm{x})} \delta_i^{-\#} \geq \eqref{Le8.ineff.strong} \ge D^{\natural -}(\bm{x},\bm{y})\\
\geq&\min\left\{  \sum_{i\in I(\bm{x})} \delta_i^{-} \,\left|  
  \begin{array}{l}
(\bm{x}-\sum_{i \in I(x)} \delta^-_ix_i\bm{e}_i ,\bm{y})\in \partial^w(T)\\
\delta_i^-=0 \ \forall i \notin I(\bm{x}), \  \bm{\delta}^-\geq \bm{0}\\ 
 \end{array} 
\right. 
\right\}
= \sum_{i\in I(\bm{x})} \delta_i^{-\#}.
\end{align*}
Therefore, \eqref{natrualDm} holds over  $ T \cap \left( \bbR^{m}_{+}\times \bbR^{s}_+\setminus\{\bm{0}\}\right)$. 
Similarly, by replacing $\partial^s(T)$ with $\partial^w(T)$ in~\eqref{LeastOutputS}, we have 
\begin{align*}
D^{\natural+}(\bm{x},\bm{y})&=\min \left\{  \,  \sum_{r\in I(\bm{y})} \delta^+_r \, \left| 
\begin{array}{l}
(\bm{x},\bm{y}+\sum_{r \in I(\bm{y})} \delta^+_r y_r \bm{e}_r )\in \partial^w(T)\\
\delta_r^+=0 \ \forall r \notin I(\bm{y}),\ \bm{\delta}^+\geq \bm{0}  
 \end{array}\right.
  \right\}.
\end{align*}
Therefore, \eqref{natrualDp} holds over  $ T \cap \left( \bbR^{m}_{+}\times \bbR^{s}_+\setminus\{\bm{0}\}\right)$.  
\par
It follows from Corollary~3 of~\cite{briec1999holder} that  
\begin{align*}
D^{\natural -}(\bm{x},\bm{y})&=\min \left\{  \sum_{i\in I(\bm{x})} \delta^-_i  \,\middle|\,  
\begin{array}{l}
(\bm{x}-\sum_{i \in I(\bm{x})} {\delta}^-_ix_i\bm{e}_i ,\bm{y})\in \partial^w(T)\\
\delta_i^-=0 \ \forall i \notin I(\bm{x}), \  \bm{\delta}^-\geq \bm{0} 
\end{array}
\right\}\\
&=\min_{i\in I(\bm{x})} \max\left\{ \delta \left|  
(\bm{x}-\delta x_i \bm{e}_i ,\bm{y})\in T 
\right. \right\} \\
\intertext{and }
D^{\natural +}(\bm{x},\bm{y})&=\min \left\{  \,  \sum_{r\in I(\bm{y})} \delta_r^+  \, \left| 
\begin{array}{l}
(\bm{x}, \bm{y}+\sum_{r \in I(\bm{y})} {\delta}^+_ry_r\bm{e}_r)\in \partial^w(T)\\
\delta_r^+=0 \ \forall r \notin I(\bm{y}), \  \bm{\delta}^+\geq \bm{0} 
\end{array}
 \right.\right\} \\
 &=
\min_{r\in I(\bm{y})} \max\left\{  \delta \left| (\bm{x},\bm{y}+\delta y_r \bm{e}_r )\in T \right.  \right\}. 
\end{align*}
\end{proof}
\subsection{Proof of~Lemma~\ref{Lemma9}}
\begin{proof}
In the proof of Theorem~\ref{Theorem4}, the following replacements are performed out:  
\begin{align*}
& T\cap \bbR^{m+s}_{++} & \to && T\cap \left( \bbR^{m}_{+}\setminus \{\bm{0}\} \times  \bbR^{s}_{+}\setminus \{\bm{0}\}\right),\\ 
&i=1,\ldots,m & \to && i \in I(\bm{x}),\\
&r=1,\ldots,s  & \to && r \in I(\bm{y}),\\
&\sum_{i=1}^m &\to && \sum_{i \in I(\bm{x})},\\
&\sum_{r=1}^s &\to && \sum_{r \in I(\bm{y})},\\
&D^{-}(\bm{x},\bm{y}) &\to && D^{\natural -}(\bm{x},\bm{y}),\\
&D^{+}(\bm{x},\bm{y}) &\to && D^{\natural +}(\bm{x},\bm{y}),\\
&F_S(\bm{x},\bm{y}) &\to && F^{\natural}_S(\bm{x},\bm{y}),\\
&\eqref{eqFs=max} &\to &&  \eqref{FsD},\\
&\eqref{LeastInput} & \to &&   \eqref{LeastInputS},\\
&\eqref{LeastOutput} & \to &&   \eqref{LeastOutputS},\\
&\eqref{BriecDm}  & \to &&  \eqref{natrualDm},\\
&\eqref{BriecDp} & \to &&  \eqref{natrualDp} .
\end{align*}
Moreover, we denote the optimal solution to the reformulated max SBM-AR model \eqref{maxSBM-ARnatural}
by $(\bm{\delta}^{-*},\bm{\delta}^{+*})$ and define
\[
(\bm{d}^{-*},\bm{d}^{+*})=\left(\sum_{i\in I(\bm{x})} \delta^{-*}_ix_i\bm{e}_i,\sum_{r\in I(\bm{y})}\delta^{+*}_ry_r\bm{e}_r\right).
\]
By replacing the proof of Theorem~\ref{Theorem4} with these replacements, we can show~\eqref{FsD}.  
Similarly, we can show~\eqref{FbD}.
\end{proof}
\subsection{Proof of~Theorem~\ref{Theorem10}}
\begin{proof}
From  \eqref{natrualDm}, \eqref{natrualDp}, and  Proposition~2 of~\cite{briec1999holder}, both  $D^{\natural-}$ and    $D^{\natural+}$ over  $T \cap \left(\bbR^{m}_{+}\setminus\{\bm{0}\}\times \bbR^{s}_+\setminus\{\bm{0}\}\right)$ are weakly monotonic inefficiency measures. 
Moreover, it follows from  \eqref{partialW=S} that $D^{\natural-}(\bm{x},\bm{y})=0$ is equivalent to $(\bm{x},\bm{y}) \in \partial^s(T)$ for $(\bm{x},\bm{y}) \in T\cap \left(\bbR^{m}_{+}\setminus\{\bm{0}\}\times \bbR^{s}_+\setminus\{\bm{0}\}\right)$.
We can show a similar property for  $D^{\natural+}$. 
Therefore, it follows from \eqref{FsD} and \eqref{FbD} that efficiency measures $F^{\natural}_S$ and  $F^{\natural}_B$ over $T \cap \left(\bbR^{m}_{+}\setminus\{\bm{0}\}\times \bbR^{s}_+\setminus\{\bm{0}\}\right)$ satisfy (I) and (WM). 
\par
In the proof of the strong monotonicity of $F_B$ in Theorem~\ref{Theorem4}, by replacing 
\begin{align*}
& T\cap \bbR^{m+s}_{++} & \to && T\cap \left( \bbR^{m}_{+}\setminus \{\bm{0}\} \times  \bbR^{s}_{+}\setminus \{\bm{0}\}\right),\\ 
&i=1,\ldots,m & \to && i \in I(\bm{x}),\\
&i'\in  \{1,\ldots,m \}& \to && i'\in I (\bm{x}),\\
&r=1,\ldots,s  & \to && r \in I(\bm{y}),\\
&r'\in \{1,\ldots,s \} & \to && r' \in I(\bm{y}),\\
&D^{-} &\to && D^{\natural -},\\
&D^{+} &\to && D^{\natural +},\\
&F_B &\to && F^{\natural}_B.
\end{align*}
we can show that efficiency measure $F^{\natural}_B$ over 
$T \cap \left(\bbR^{m}_{+}\setminus\{\bm{0}\}\times \bbR^{s}_+\setminus\{\bm{0}\}\right)$
satisfies (SM) if $m\leq s$.
\end{proof}
\subsection{Proof of~Theorem~\ref{Theorem11}}
\begin{proof}
Suppose $(\bm{x},\bm{y})\in T \cap \left(\bbR^{m}_{+}\setminus\{\bm{0}\}\times \bbR^{s}_+\setminus\{\bm{0}\}\right)$, then 
there exists an index $\bar{i}\in I(\bm{x})$ such that 
\begin{align}
D^{\natural -}(\bm{x},\bm{y})=\max\left\{ \delta \left| (\bm{x}-\delta x_{\bar{i}} \bm{e}_{\bar{i}},\bm{y})\in T\right\}.\right.\label{eqibar}
\end{align}
The problem $\max\left\{ \delta \left| (\bm{x}-\delta x_{\bar{i}} \bm{e}_{\bar{i}},\bm{y})\in T\right\}\right.$
 is a dual problem of  
\[
\min\left\{ \bm{v}\bm{x}-\bm{u}\bm{y} \left|
\begin{array}{l} 
\bm{v}\bm{x}_j - \bm{u}\bm{y}_j \geq 0 \ (j=1,\ldots,n), \\ 
-\bm{v}P \geq \bm{0}, -\bm{u}Q \geq \bm{0},\ v_{\bar{i}}x_{\bar{i}}=1,\\
\bm{v}\ge\bm{0},\ \bm{u}\ge\bm{0}
\end{array}
\right\}.\right.
\]
Suppose that for any $i\in I(\bm{x})$, there are $L_i$ vertices,
$(\bm{v}^1_i,\bm{u}^1_i)\cdots,(\bm{v}^{L_i}_i,\bm{u}^{L_i}_i)$,  for each polyhedron, 
\[
\left\{ (\bm{v},\bm{u}) \left| 
\begin{array}{l}
\bm{v}\bm{x}_j - \bm{u}\bm{y}_j \geq 0 \ (j=1,\ldots,n),\\
-\bm{v}P \geq \bm{0}, -\bm{u}Q \geq \bm{0},\ v_ix_i=1,\\
\bm{v}\ge\bm{0},\ \bm{u}\ge\bm{0}
\end{array}
\right\}.\right.
\] 
Without loss of generality, assume that  $D^{\natural -}(\bm{x},\bm{y})=\bm{v}_{\bar{i}}^1\bm{x}-\bm{u}_{\bar{i}}^1\bm{y}$.
\par
Suppose any $(\bm{\epsilon}^-,\bm{\epsilon}^+)\in \bbR^{m+s}_{+}$ such that $\epsilon_i^->0 \,(i\notin I(\bm{x}))$,  $\epsilon_i^-=0 \,(i\in I(\bm{x}))$,  $\epsilon_r^+>0 \,(r\notin I(\bm{y}))$ and 
$\epsilon_r^+=0 \,(r\in I(\bm{y}))$, then $\bm{x}(\bm{\epsilon}^-)=\bm{x}+\bm{\epsilon}^-$,  $\bm{y}(\bm{\epsilon}^+)=\bm{y}+\bm{\epsilon}^+$ and   
\begin{align}
 D^-(\bm{x}(\bm{\epsilon}^-),\bm{y}(\bm{\epsilon}^+))&= D^{\natural-}(\bm{x}(\bm{\epsilon}^-),\bm{y}(\bm{\epsilon}^+))\\
 & \leq   \max\left\{ \delta \left| (\bm{x}(\bm{\epsilon}^-)-\delta x_{\bar{i}} \bm{e}_{\bar{i}},\bm{y}(\bm{\epsilon}^+))\in T\right\}\right.\\
 &= \min\left\{ \bm{v}\bm{x}(\bm{\epsilon}^-)-\bm{u}\bm{y}(\bm{\epsilon}^+) \left|
\begin{array}{l} 
\bm{v}\bm{x}_j - \bm{u}\bm{y}_j \geq 0 \ (j=1,\ldots,n) \\ 
-\bm{v}P \geq \bm{0}, -\bm{u}Q \geq \bm{0},\\ v_{\bar{i}}x_{\bar{i}}=1 ,
\bm{v}\ge\bm{0},\ \bm{u}\ge\bm{0}
\end{array}
\right\}.\right.\\
& \leq \bm{v}_{\bar{i}}^1\bm{x}(\bm{\epsilon}^-)-\bm{u}_{\bar{i}}^1\bm{y}(\bm{\epsilon}^+)  \leq \bm{v}_{\bar{i}}^1\bm{x}+ \bm{v}_{\bar{i}}^1\bm{\epsilon}^- -\bm{u}_{\bar{i}}^1\bm{y}\\
&= D^{\natural -}(\bm{x},\bm{y})+ \bm{v}_{\bar{i}}^1\bm{\epsilon}^- \label{eqLeftW}
\intertext{and}
 D^-(\bm{x}(\bm{\epsilon}^-),\bm{y}(\bm{\epsilon}^+))&= D^{\natural-}(\bm{x}(\bm{\epsilon}^-),\bm{y}(\bm{\epsilon}^+)) \geq  D^{\natural-}(\bm{x},\bm{y}(\bm{\epsilon}^+))\\
 &=\min\left\{ \bm{v}^l_k\bm{x}-\bm{u}^l_k\bm{y}(\bm{\epsilon}^+) \left|\begin{array}{l}  k\in I(\bm{x}) \\ l \in \{1,\ldots, L_k \} \end{array}  \right. \right\}\\
 &= \min\left\{ \bm{v}^l_k\bm{x}-\bm{u}^l_k\bm{y}- \bm{u}^l_k\bm{\epsilon}^+ \left|
 \begin{array}{l}  k\in I(\bm{x}) \\ l \in \{1,\ldots, L_k \} \end{array}  \right. \right\}\\
 &\geq  \min\left\{ \bm{v}^l_k\bm{x}-\bm{u}^l_k\bm{y} \left| 
 \begin{array}{l}  k\in I(\bm{x}) \\ l \in \{1,\ldots, L_k \} \end{array} 
\right. \right\} \\
&\qquad \qquad -\max\left\{ \bm{u}^l_k\bm{\epsilon}^+ \left| 
\begin{array}{l}  k\in I(\bm{x}) \\ l \in \{1,\ldots, L_k \} \end{array} \right. \right\} \\
&=D^{\natural -}(\bm{x},\bm{y})-\max \left\{ \bm{u}^l_k\bm{\epsilon}^+ \left| 
\begin{array}{l}  k\in I(\bm{x}) \\ l \in \{1,\ldots, L_k \} \end{array} \right. \right\}. \label{eqRightW}
 \end{align}
 Since 
 \[
 \lim_{\bm{\epsilon}^-\to +\bm{0} \atop \epsilon^-_i=0\ i \in I(\bm{x})}  \bm{v}_{\bar{i}}^1\bm{\epsilon}^-=
 \lim_{\bm{\epsilon}^+\to +\bm{0} \atop  \epsilon_r^+=0\ r\in I(\bm{y})}  \max \left\{ \bm{u}^l_k\bm{\epsilon}^+ \left| 
\begin{array}{l}  k\in I(\bm{x}) \\ l \in \{1,\ldots, L_k \} \end{array} \right. \right\}=0,  
 \]
 it follows from  \eqref{eqLeftW}  and \eqref{eqRightW} that 
\[
D^{\natural-}(\bm{x},\bm{y}) =  \lim_{\mathscr{E}(\bm{x},\bm{y})\ni (\bm{\epsilon}^-,\bm{\epsilon}^+) \to (\bm{0},\bm{0}) } D^-(\bm{x}(\bm{\epsilon}^-),\bm{y}(\bm{\epsilon}^+)).
\]
Similarly, we can show that
\[
D^{\natural+}(\bm{x},\bm{y}) =  \lim_{ \mathscr{E}(\bm{x},\bm{y})\ni (\bm{\epsilon}^-,\bm{\epsilon}^+)\to (\bm{0},\bm{0}) } D^+(\bm{x}(\bm{\epsilon}^-),\bm{y}(\bm{\epsilon}^+)). 
\]
\par
From~\eqref{FsD},~\eqref{FbD} of Lemma~\ref{Lemma9}, and $D^{\natural +}(\bm{x},\bm{y})\geq 0$, functions $F_B^{\natural}$ and $F_S^{\natural}$ are 
continuous over $T \cap \left(\bbR^{m}_{+}\setminus\{\bm{0}\}\times \bbR^{s}_+\setminus\{\bm{0}\}\right)$; hence, \eqref{FSnat} and \eqref{FBnat} hold.
\end{proof}

\end{document}